\documentclass{article}
\usepackage{amsmath,amsthm}
\usepackage{verbatim,amsmath,amscd,amssymb,amsthm}
\begin{document}
\font\germ=eufm10
\def\ssl{\hbox{\germ sl}}
\def\slh{\widehat{\ssl_2}}
\def\ge{\hbox{\germ g}}
\def\aaa{@}
\title{\Large\bf Evaluation representations 
of quantum affine algebras \\
at roots of unity }

\author{
Yuuki A\textsc{be} 
\thanks
{
e-mail: yu-abe@sophia.ac.jp
}
\\
Department of Mathematics, 
\\
Sophia University
\and
Toshiki N\textsc{akashima} 
\thanks
{
supported in part by JSPS Grants in Aid for 
Scientific Research, 
e-mail: toshiki@mm.sophia.ac.jp
}
\\
Department of Mathematics, 
\\
Sophia University
}
\date{}
\maketitle
\begin{abstract}
The purpose of this paper is 
to compute the Drinfel'd polynomials
for two types of evaluation representations 
of quantum affine algebras at roots of unity 
and construct those representations 
as the submodules of 
evaluation Schnizer modules. 
Moreover, we obtain
the necessary and sufficient condition for that
the two types of evaluation representations 
are isomorphic to each other.
\end{abstract}
\maketitle

\renewcommand{\labelenumi}{$($\roman{enumi}$)$}
\renewcommand{\labelenumii}{$(${\rm \alph{enumii}}$)$}
\font\germ=eufm10
\newcommand{\cA}{{\mathcal A}}
\newcommand{\cB}{{\mathcal B}}
\newcommand{\cC}{{\mathcal C}}
\newcommand{\cD}{{\mathcal D}}
\newcommand{\cF}{{\mathcal F}}
\newcommand{\cH}{{\mathcal H}}
\newcommand{\cI}{{\mathcal I}}
\newcommand{\cK}{{\mathcal K}}
\newcommand{\cL}{{\mathcal L}}
\newcommand{\cM}{{\mathcal M}}
\newcommand{\cN}{{\mathcal N}}
\newcommand{\cO}{{\mathcal O}}
\newcommand{\cS}{{\mathcal S}}
\newcommand{\cV}{{\mathcal V}}
\newcommand{\cW}{{\mathcal W}}
\newcommand{\fra}{\mathfrak a}
\newcommand{\frb}{\mathfrak b}
\newcommand{\frc}{\mathfrak c}
\newcommand{\frd}{\mathfrak d}
\newcommand{\fre}{\mathfrak e}
\newcommand{\frf}{\mathfrak f}
\newcommand{\frg}{\mathfrak g}
\newcommand{\frh}{\mathfrak h}
\newcommand{\fri}{\mathfrak i}
\newcommand{\frj}{\mathfrak j}
\newcommand{\frk}{\mathfrak k}
\newcommand{\fm}{\mathfrak m}
\newcommand{\frn}{\mathfrak n}
\newcommand{\frp}{\mathfrak p}
\newcommand{\fq}{\mathfrak q}
\newcommand{\frr}{\mathfrak r}
\newcommand{\frs}{\mathfrak s}
\newcommand{\frt}{\mathfrak t}
\newcommand{\fru}{\mathfrak u}
\newcommand{\frA}{\mathfrak A}
\newcommand{\frB}{\mathfrak B}
\newcommand{\frF}{\mathfrak F}
\newcommand{\frG}{\mathfrak G}
\newcommand{\frH}{\mathfrak H}
\newcommand{\frI}{\mathfrak I}
\newcommand{\frJ}{\mathfrak J}
\newcommand{\frN}{\mathfrak N}
\newcommand{\frP}{\mathfrak P}
\newcommand{\frT}{\mathfrak T}
\newcommand{\frU}{\mathfrak U}
\newcommand{\frV}{\mathfrak V}
\newcommand{\frX}{\mathfrak X}
\newcommand{\frY}{\mathfrak Y}
\newcommand{\frZ}{\mathfrak Z}
\newcommand{\rA}{\mathrm{A}}
\newcommand{\rC}{\mathrm{C}}
\newcommand{\rd}{\mathrm{d}}
\newcommand{\rB}{\mathrm{B}}
\newcommand{\rD}{\mathrm{D}}
\newcommand{\rE}{\mathrm{E}}
\newcommand{\rH}{\mathrm{H}}
\newcommand{\rK}{\mathrm{K}}
\newcommand{\rL}{\mathrm{L}}
\newcommand{\rM}{\mathrm{M}}
\newcommand{\rN}{\mathrm{N}}
\newcommand{\rR}{\mathrm{R}}
\newcommand{\rT}{\mathrm{T}}
\newcommand{\rZ}{\mathrm{Z}}
\def\al{\alpha}
\def\b{\beta}
\def\d{\delta}
\def\D{\Delta}
\def\e{\varepsilon}
\def\ep{\epsilon}
\def\ga{\gamma}
\def\l{\lambda}
\def\L{\Lambda}
\def\Om{\Omega}
\def\t{\theta}
\def\vep{\varepsilon}
\def\vp{\varphi}
\newcommand{\bbA}{\mathbb A}
\newcommand{\bbC}{\mathbb C}
\newcommand{\bbG}{\mathbb G}
\newcommand{\bbF}{\mathbb F}
\newcommand{\bbH}{\mathbb H}
\newcommand{\bbP}{\mathbb P}
\newcommand{\bbN}{\mathbb N}
\newcommand{\bbQ}{\mathbb Q}
\newcommand{\bbR}{\mathbb R}
\newcommand{\bbV}{\mathbb V}
\newcommand{\bbZ}{\mathbb Z}
\newcommand{\adj}{\operatorname{adj}}
\newcommand{\Ad}{\mathrm{Ad}}
\newcommand{\Ann}{\mathrm{Ann}}
\newcommand{\rcris}{\mathrm{cris}}
\newcommand{\ch}{\mathrm{ch}}
\newcommand{\coker}{\mathrm{coker}}
\newcommand{\diag}{\mathrm{diag}}
\newcommand{\Diff}{\mathrm{Diff}}
\newcommand{\Dist}{\mathrm{Dist}}
\newcommand{\rDR}{\mathrm{DR}}
\newcommand{\ev}{\mathrm{ev}}
\newcommand{\Ext}{\mathrm{Ext}}
\newcommand{\cExt}{\mathcal{E}xt}
\newcommand{\fin}{\mathrm{fin}}
\newcommand{\Frac}{\mathrm{Frac}}
\newcommand{\GL}{\mathrm{GL}}
\newcommand{\Hom}{\mathrm{Hom}}
\newcommand{\hd}{\mathrm{hd}}
\newcommand{\rht}{\mathrm{ht}}
\newcommand{\id}{\mathrm{id}}
\newcommand{\im}{\mathrm{im}}
\newcommand{\inc}{\mathrm{inc}}
\newcommand{\ind}{\mathrm{ind}}
\newcommand{\coind}{\mathrm{coind}}
\newcommand{\Lie}{\mathrm{Lie}}
\newcommand{\Max}{\mathrm{Max}}
\newcommand{\mult}{\mathrm{mult}}
\newcommand{\op}{\mathrm{op}}
\newcommand{\ord}{\mathrm{ord}}
\newcommand{\pt}{\mathrm{pt}}
\newcommand{\qt}{\mathrm{qt}}
\newcommand{\rad}{\mathrm{rad}}
\newcommand{\res}{\mathrm{res}}
\newcommand{\rgt}{\mathrm{rgt}}
\newcommand{\rk}{\mathrm{rk}}
\newcommand{\SL}{\mathrm{SL}}
\newcommand{\soc}{\mathrm{soc}}
\newcommand{\Spec}{\mathrm{Spec}}
\newcommand{\St}{\mathrm{St}}
\newcommand{\supp}{\mathrm{supp}}
\newcommand{\Tor}{\mathrm{Tor}}
\newcommand{\Tr}{\mathrm{Tr}}
\newcommand{\wt}{\mathrm{wt}}
\newcommand{\Ab}{\mathbf{Ab}}
\newcommand{\Alg}{\mathbf{Alg}}
\newcommand{\Grp}{\mathbf{Grp}}
\newcommand{\Mod}{\mathbf{Mod}}
\newcommand{\Sch}{\mathbf{Sch}}\newcommand{\bfmod}{{\bf mod}}
\newcommand{\Qc}{\mathbf{Qc}}
\newcommand{\Rng}{\mathbf{Rng}}
\newcommand{\Top}{\mathbf{Top}}
\newcommand{\Var}{\mathbf{Var}}
\newcommand{\gromega}{\langle\omega\rangle}
\newcommand{\lbr}{\begin{bmatrix}}
\newcommand{\rbr}{\end{bmatrix}}
\newcommand{\forb}{\bigcirc\kern-2.8ex \because}
\newcommand{\forbb}{\bigcirc\kern-3.0ex \because}
\newcommand{\forbbb}{\bigcirc\kern-3.1ex \because}
\newcommand{\SpS}{spectral sequence}
\newcommand\C{\mathbb C}
\newcommand\hh{{\hat{H}}}
\newcommand\eh{{\hat{E}}}
\newcommand\F{\mathbb F}
\newcommand\fh{{\hat{F}}}
\newcommand{\End}{\operatorname{End}}
\newcommand{\Stab}{\operatorname{Stab}}
\newcommand{\mo}{\operatorname{mod}}
\newcommand\pf{\noindent {\bf Proof:  }}
\def\AA{{\cal A}}
\def\arr{\longrightarrow}
\def\bf{\textbf}
\def\Bfin{B_{\e}^{\rm{fin}}}
\def\DD{{\cal D}}
\def\Di{\tilde{\D}^{\textrm{im}}}
\def\Dr{\tilde{\D}^{\textrm{re}}}
\def\g{\frg}
\def\gl{{\mathfrak g}{\mathfrak l}}
\def\h{\frh}
\def\hIfin{\wh{I}_{\e}^{\rm{fin}}}
\def\hMfin{\wh{M}_{\e}^{\rm{fin}}}
\def\hNfin{\wh{N}_{\e}^{\rm{fin}}}
\def\hUfin{\wh{U}_{\e}^{\rm{fin}}}
\def\it{\textit}
\def\Lfin{L_{\e}^{\textrm{fin}}(\l)}
\def\Lnil{L_{\e}^{\textrm{nil}}(\l)}
\def\nil{L_{l}^{nil}(\l)}
\def\no{\nonumber}
\def\ot{\otimes}
\def\PP {{\mathcal P}}
\def\q{\quad}
\def\qq{\qquad}
\def\rm{\textrm}
\def\RR {{\mathcal R}}
\def\s={\star \in \{-,0,+\}}
\def\sl{{\mathfrak s}{\mathfrak l}}
\def\SS{{\cal S}}
\def\ta{\widetilde{a}}
\def\tb{\widetilde{b}}
\def\tB{\widetilde{B}}
\def\tBfin{\tB_{\e}^{\rm{fin}}}
\def\tD{\widetilde{\D}}
\def\tep{\widetilde{\epsilon}}
\def\tg{\widetilde{g}}
\def\tI{\widetilde{I}}
\def\tm{\widetilde{m}}
\def\tMfin{\wt{M}_{\e}^{\rm{fin}}}
\def\tNfin{\wt{N}_{\e}^{\rm{fin}}}
\def\TT{{\cal T}}
\def\tIfin{\wt{I}_{\e}^{\rm{fin}}}
\def\tu{\widetilde{U}}
\def\tU{\widetilde{U}}
\def\tue{\widetilde{U}_{\e}}
\def\tuf{\widetilde{U}_{\e}^{\rm{fin}}}
\def\tUe{\wt{U}_{\e}}
\def\tUfin{\widetilde{U}_{\e}^{\rm{fin}}}
\def\tUres{\widetilde{U}_{\e}^{\rm{res}}}
\def\tV{\widetilde{V}}
\def\tv{\widetilde{v}}
\def\tVdia{\wt{V}^{\rm{dia}}}
\def\tVfin{\tV^{\textrm{fin}}_{\e}(\bf{P})}
\def\tVnil{\tV^{\textrm{nil}}_{\e}(\bf{P})}
\def\tVres{\tV^{\textrm{res}}_{\e}(\bf{P})}
\def\tw{\widetilde{W}}
\def\tx{\widetilde{x}}
\def\ty{\widetilde{y}}
\def\tz{\widetilde{z}}
\def\ua{U_{\varepsilon}(A_n)}
\def\ub{U_{\varepsilon}(B_n)}
\def\uc{U_{\varepsilon}(C_n)}
\def\ud{U_{\varepsilon}(D_n)}
\def\ue{U_{\varepsilon}}
\def\Ue{U_{\e}}
\def\uf{U_{\varepsilon} ^{\rm{fin}}}
\def\Ufin{U_{\varepsilon} ^{\rm{fin}}}
\def\uq{U_q(\mathfrak g)}
\def\ur{U_{\varepsilon} ^{\rm{res}}}
\def\Ures{U_{\varepsilon} ^{\rm{res}}}
\def\Vfin{V^{\textrm{fin}}_{\e}(\l)}
\def\Vnil{V^{\textrm{nil}}_{\e}(\l)}
\def\Vres{V^{\textrm{res}}_{\e}(\l)}
\def\wh{\widehat}
\def\wt{\widetilde}
\def\WW {{\mathcal W}}
\section{Introduction} 

For a generic $q$, 
let $U_q(\g)$ be the quantum algebra associated with 
a simple Lie algebra $\g$ 
and $U_q(\wt{\g})$ be the non-twisted quantum loop algebra of $\g$. 
It is known that every finite dimensional irreducible $U_q(\g)$ 
(resp. $U_q(\wt{\g})$) modules are highest weight module 
and classified by highest weights. 
Moreover, there exists one to one correspondence from 
the set of their highest weights to $\bbZ_{+}^n$ 
(resp. polynomials $\bbC_0[t]^n$), 
where $\bbZ_{+}:=\{0,1,2, \cdots\}$ 
(resp. $\bbC_0[t]:=\{P \in \bbC[t] 
\, | \, P \rm{ is monic and } P(0) \neq 0\}$). 
The theory of finite dimensional $U_q(\wt{\g})$-modules 
is introduced in \cite{CP95}. 
We denote the $U_q(\g)$ (resp. $U_q(\wt{\g})$) module corresponding to 
$\l \in \bbZ_{+}^n$ (resp. $\bf{P} \in \bbC_0[t]^n$) by 
$V_q(\l)$ (resp. $\tV_q(\bf{P})$), 
where the polynomial $\bf{P}$ of $\tV_q(\bf{P})$ is called 
``Drinfel'd polynomial''. 

In the case $\g=\sl_{n+1}$, 
there exist $\bbC$-algebra homomorphisms 
$\rm{ev}_{\bf{a}}^{+}, \rm{ev}_{\bf{a}}^{-}
: U_q(\wt{\sl}_{n+1}) \arr U_q(\sl_{n+1})$ 
for $\bf{a} \in \bbC^{\times}$ (see \cite{J}, \cite{CP94a}). 
By using these homomorphisms, 
we can regard $V_q(\l)$ as a $U_q(\wt{\sl}_{n+1})$-module, 
which are called ``evaluation representations'' 
and denoted by $V_q(\l)_{\bf{a}}^{\pm}$. 
By the classification theorem of finite dimensional 
$U_q(\wt{\sl}_{n+1})$-modules(\cite{CP95}), 
there exists a unique polynomial
$\bf{P}_{\bf{a}}^{\pm} \in \bbC_0[t]^n$ such that 
$V_q(\l)_{\bf{a}}^{\pm}$ is isomorphic to $\tV_q(\bf{P}_{\bf{a}}^{\pm})$ 
as a $U_q(\wt{\sl}_{n+1})$-module respectively.
The Drinfel'd polynomials $\bf{P}_{\bf{a}}^{\pm}$ are 
computed by Chari and Pressley in \cite{CP94a}. 
In this paper, we shall consider evaluation representations 
at roots of unity.  

Let $\e$ be a primitive $l$-th root of unity. 
The representation theory of quantum algebras at roots of unity 
is divided into two types. 
One is for $U_{\e}(\g)$, $U_{\e}(\wt{\g})$ defined by 
De Concini-Kac (=non-restricted type) in \cite{DK} 
and the other is for $\Ures(\g)$, $\Ures(\wt{\g})$ defined by 
Lusztig (=restricted type) in \cite{L89}. 

$\Ures(\g)$ (resp. $\Ures(\wt{\g})$) 
has the $\bbC$-subalgebra $\Ufin(\g)$ (resp. $\Ufin(\wt{\g})$) 
which is called ``small quantum algebra''. 
By the tensor product theorem (see \cite{L89}, \cite{CP97}), 
in order to understand the representation theory of 
$\Ures(\g)$ (resp. $\Ures(\wt{\g})$), 
we may consider the one of $\Ufin(\g)$ (resp. $\Ufin(\wt{\g})$). 
Indeed, every finite dimensional irreducible $\Ufin(\g)$ 
(resp. $\Ufin(\wt{\g})$) module is a highest weight module 
and classified by highest weight.
Moreover, there exists one to one correspondence from the set
of their highest weights to $\bbZ_l^n$ 
(resp. polynomials $\bbC_l[t]^n$), where $\bbZ_l:=\{0,1, \cdots l-1\}$ 
(resp. $\bbC_l[t]:=\{P \in \bbC_0[t] \, | \, P \rm{ is not divisible by } 
(1-ct^l) \rm{ for all } c \in \bbC^{\times}\}$). 
We denote the $\Ufin(\g)$ (resp. $\Ufin(\wt{\g})$) module 
corresponding to $\l \in \bbZ_{l}^n$ (resp. $\bf{P} \in \bbC_l[t]^n$) by 
$\Vfin$ (resp. $\tVfin$). \\
\q We also obtain the evaluation representations of $\Vfin$
in the case of $\Ufin(\wt{\sl}_{n+1})$. 
We denote them by $\Vfin_{\bf{a}}^{\pm}$. 
We can compute the Drinfel'd polynomials of $\Vfin_{\bf{a}}^{\pm}$ 
by the similar method to \cite{CP94a} 
(see Theorem 4.13 in this paper). 
Moreover, for $\bf{a}_{\pm} \in \bbC^{\times}$, we shall show
 that 
$\Vfin_{\bf{a}_{+}}^{+}$ is isomorphic to $\Vfin_{{\bf{a}}_{-}}^{-}$ 
if and only if 
\begin{eqnarray}
{\bf{a}}_{+}=\bf{a}_{-}\e^{2(\sum_{k=1}^{i-1}\l_k-\sum_{k=i+1}^n\l_k+i)} 
\rm{  for all }i \in \rm{supp}(\l),
\label{intro}
\end{eqnarray}
where $\rm{supp}(\l):=\{1 \leq i \leq n\, |\l_i \neq 0\}$.
If $q$ is generic, the condition (1) never occurs for $\#(\rm{supp}(\l))>1$. 
But, in this case, there exists  $\l \in \bbZ_l^n$  
which satisfies (\ref{intro}) for $\#(\rm{supp}(\l))>1$ 
(see Proposition \ref{pro pc=pc22}, \ref{pro pc=pc222}).  

On the other hand, many finite-dimensional irreducible 
$U_{\e}(\g)$ (resp. $U_{\e}(\wt{\g})$) modules 
are no longer highest or lowest weight modules 
and they are characterized by several continuous parameters 
(see \cite{DK}, \cite{BK}). 
For $\g=\sl_{n+1}$, such $U_{\e}(\sl_{n+1})$-modules 
are constructed explicitly in \cite{DJMM}, 
which are called ``maximal cyclic representation''. 
For an arbitrary simple Lie algebra $\g$, 
Schnizer introduced an alternative construction of 
such $U_{\e}(\g)$-modules in \cite{S93}, \cite{S94}, 
which we call ``Schnizer modules''. 

By using the theory of 
the quantum algebra of restricted type, 
we obtain that every finite dimensional 
irreducible ``nilpotent'' $U_{\e}(\sl_{n+1})$-modules are 
highest weight module and classified by highest weight 
(see \S3.4 and \S 5.2). 
Moreover, there exists the one to one
correspondence from the set of their highest weights to 
$\bbZ_l^n$ (resp. $\bbC_l[t]^n$).
We denote the $U_{\e}(\sl_{n+1})$ (resp. $U_{\e}(\wt{\sl}_{n+1})$) module 
corresponding to $\l \in \bbZ_{l}^n$
(resp. $\bf{P} \in \bbC_l[t]^n$) by $V^{\rm{nil}}_{\e}(\l)$
(resp. $\wt{V}^{\rm{nil}}_{\e}(\bf{P})$). 
We also obtain the evaluation representations of $\Vnil$,
which are denoted by $\Vnil_{\bf{a}}^{\pm}$. 
The module $\Vnil_{\bf{a}}^{\pm}$ is regraded as a 
$\Ufin(\wt{\sl}_{n+1})$-module 
and $\Vnil_{\bf{a}}^{\pm}$ is isomorphic to 
$\Vfin_{\bf{a}}^{\pm}$ as a 
$\Ufin(\wt{\sl}_{n+1})$-module (see \S 5.2). 
Therefore, for $\bf{a}_{\pm} \in \bbC^{\times}$, we obtain  
\begin{eqnarray}
\rm{$\Vnil_{\bf{a}_{+}}^{+}$ is isomorphic to $\Vnil_{{\bf{a}}_{-}}^{-}$ 
if and only if (1) holds}. 
\label{intro2}
\end{eqnarray} 
\q We can also prove (\ref{intro2}) without using 
the theory of the quantum algebra of restricted type. 
In \cite{N}, T.N. showed that one can 
construct $\Vnil$ as the subrepresentation of a maximal cyclic
representation by specializing their parameters properly 
for type $A$.
Similarly, in \cite{AN}, we found that we can construct $\Vnil$ 
as a submodule of a Schnizer module if $\g$=A, B, C or D, 
and then we can construct $\Vnil_{\bf{a}}^{\pm}$ as 
the submodule of evaluation of a Schnizer module. 
By using this fact, we can prove (\ref{intro2})
(see \S 5 alternative proof of Proposition 
\ref{pro isomorphic condition} (b)). 

The organization of this paper is as follows. 
In \S 2, we introduce basic properties of quantum algebras for generic $q$. 
In \S 3, we introduce quantum algebras at roots of unity 
of non-restricted type and restricted type. 
Moreover, we prove the isomorphisn theorem of these algebras. 
In \S4 (resp. \S 5), we discuss about the evaluation representations 
of restricted (resp. non-restricted) type. 
\section{Quantum algebras (generic case)}
\setcounter{equation}{0}
\renewcommand{\theequation}{\thesection.\arabic{equation}}
\subsection{Notations}
\q We fix the following notations (see \cite{B94a}, \cite{BK}). 
Let $\sl_{n+1}$ be the finite dimensional simple Lie algebra 
over $\bbC$ of type $A_n$ 
and $\widetilde{\sl}_{n+1}=\sl_{n+1} \otimes \bbC[t,t^{-1}]$ 
the loop algebra of $\sl_{n+1}$. 
We set $I:=\{1,2, \cdots , n\}$ and $\widetilde{I}:=I \sqcup \{0\}$.
Let $(\fra_{i,j})_{i,j \in \tI}$ 
be the generalized Cartan matrix of $\widetilde{\sl}_{n+1}$, 
that is, $\fra_{i,i}=2$, $\fra_{i,j}=-1$ if $|i-j|=1$ or $n$, 
and $\fra_{i,j}=0$ otherwise. 
Then $(\fra_{i,j})_{i,j \in I}$ 
is the Cartan matrix of $\sl_{n+1}$. 
Let $\Pi:=\{\al_i\}_{i \in I}$ (resp. $\wt{\Pi}:=\{\al_i\}_{i \in \tI}$) 
the set of the simple roots of $\sl_{n+1}$ (resp. $\wt{\sl}_{n+1}$) 
and $\Pi^{\vee}:=\{\al_i^{\vee}\}_{i \in I}$ 
(resp. $\wt{\Pi}^{\vee}:=\{\al_i^{\vee}\}_{i \in \tI}$) 
be the set of the simple coroots of $\sl_{n+1}$ (resp. $\wt{\sl}_{n+1}$). 
Let $\h$ be the Cartan subalgebra of $\sl_{n+1}$ 
and $\h^{*}$ the $\bbC$-dual space of $\h$. 
Then $\Pi^{\vee}$ (resp. $\Pi$ ) is a $\bbC$-basis of $\h$
(resp. $\h^{*}$). 
We have a $\bbC$-bilinear map 
$\langle, \rangle: \h^{*} \times \h \longrightarrow \bbC$ 
such that $\langle \al_j, \al_i^{\vee} \rangle =\fra_{i,j}$ 
for any $i,j \in I$. 
Define the root lattice $Q:=\bigoplus_{i \in I}\bbZ \al_i$ 
(resp. the coroot lattice 
$Q^{\vee}:=\bigoplus_{i \in I}\bbZ \al_i^{\vee}$) 
and the affine root lattice $\wt{Q}:=\bbZ\al_0\oplus Q$ 
(resp. the affine coroot lattice 
$\wt{Q}^{\vee}:=\bbZ\al_0^{\vee} \oplus Q^{\vee}$).
For $i \in I$, we define the fundamental weights 
$\{\L_i\}_{i \in I} \subset \h^{*}$ by
\begin{eqnarray}
\L_i:=\frac{1}{n+1}
\{(n-i+1)\sum_{k=1}^i k\al_k+i\sum_{k=i+1}^n(n-k+1)\al_k\}. 
\label{def fundamental weights}
\end{eqnarray}
Similarly, we define the fundamental coweights 
$\{\L_i^{\vee}\}_{i \in I} \subset \h$ 
by replacing $\al$ in $\L_i$ with $\al^{\vee}$. 
Then we have 
$\langle \L_i, \al_j^{\vee} \rangle=\delta_{i,j}$
(resp. $\langle \al_j, \L_i^{\vee} \rangle=\delta_{i,j}$) 
for any $i,j \in I$. 
Define the weight lattice 
$P:=\bigoplus_{i \in I}\bbZ\L_i$ 
(resp. the coweight lattice $P^{\vee}
:=\bigoplus_{i \in I}\bbZ\L_i^{\vee}$)
and define a symmetric bilinear form 
$(,): \h^{*} \times \h^{*} \longrightarrow \bbC$ 
determined by $(\al_i, \al_j)= \fra_{i,j}$ for any $i,j \in I$. \\
\q Let $\D$ (resp. $\D_{+}$) be the set of roots 
(resp. positive roots) of $\sl_{n+1}$ 
and $\t:=\sum_{i \in I}\al_i$ be the highest root in $\D$. 
We set $\delta:=\al_0+\t$. 
Let $\widetilde{\D}$ be the affine root system of $\wt{\sl}_{n+1}$. 
Then we have $\widetilde{\D}=\widetilde{\D}^{\textrm{re}}
 \sqcup \widetilde{\D}^{\textrm{im}}$, 
where 
\begin{eqnarray*}
\widetilde{\D}^{\textrm{re}}:=\{\al+n\delta \, | \, 
\al \in \D , n \in \bbZ\},
\q \widetilde{\D}^{\textrm{im}}:=\{n\delta \, | \,  n \in 
\bbZ^{\times}:=(\bbZ\setminus\{0\})\}, 
\end{eqnarray*}
and  $\tD=\tD_{+} \sqcup (-\tD_{+})$, where 
\begin{eqnarray*}
\widetilde{\D}_{+}:=\widetilde{\D}_{+}^{\textrm{re}}
 \sqcup \widetilde{\D}_{+}^{\textrm{im}}, 
\q \widetilde{\D}_{+}^{\textrm{re}}:=\{\al+n\delta \, | \, 
\al \in \D , n \in \bbN:=\{1,2, \cdots\}\} \sqcup \D_{+},
\q \widetilde{\D}_{+}^{\textrm{im}}:=\{n\delta \, | \,  n \in \bbN\}. 
\end{eqnarray*} 
Moreover, we set 
\begin{eqnarray*}
&&\tD_{+}^{\textrm{im}}(I):= I \times \tD_{+}^{\textrm{im}}
=\{(i, n \delta) \, | \, i \in I, n \in \bbN\}, 
\q \tD_{+}(I):=\tD_{+}^{\textrm{re}} \sqcup \tD_{+}^{\textrm{im}}(I), \\
&&\tD^{\textrm{im}}(I):=\{(i, n \delta) \, | \, i \in I, n \in
 \bbZ^{\times}\}, 
\q \tD(I):=\tD^{\textrm{re}} \sqcup \tD^{\textrm{im}}(I).
\end{eqnarray*} 
\q For $i \in \widetilde{I}$, 
let $s_{i}$ be the simple reflection on $\h^{*}$, 
that is,  $s_{i}(\l)=\l- \langle \l, \al_i^{\vee} 
\rangle\al_i$ for any $\l \in \h^{*}$. 
The affine Weyl group $\widetilde{\WW}$ of $\widetilde{\sl}_{n+1}$ 
(resp. Weyl group $\WW$ of $\sl_{n+1}$) 
is generated by $\{s_i\}_{i \in \widetilde{I}}$ 
(resp. $\{s_i\}_{i \in I}$). 
For $x \in \h$, 
we define $t_x:\h^{*} \longrightarrow \h^{*}$ 
by $t_x(\l)=\l-\langle \l, x \rangle \delta$ 
and set $T_{P^{\vee}}:=\{t_x \, | \, x \in P^{\vee}\}$, 
$T_{Q^{\vee}}:=\{t_x \, | \, x \in Q^{\vee}\}$.
Consider the extended affine Weyl group 
$\widehat{\WW}:= \WW  \widetilde{\times} T_{P^{\vee}} $, 
where the structure of the semi-direct product is given by 
$(s,t_x)(s^{'},t_y)=(ss^{'}, t_{s^{'-1}x}t_y)$ 
for any $s,s^{'} \in \WW$, and $x,y \in P^{\vee}$. 
We set $\TT:=\{\tau: \widetilde{I} \arr \tI; \textrm{ permutation } 
| \, \fra_{\tau(i), \tau(j)}=\fra_{i,j} \,\textrm{for any}\, i,j \in
\tI\}$ 
and define the semi-direct product 
$\TT \widetilde{\times} \widetilde{\WW}$ 
by $\tau s_i \tau^{-1}=s_{\tau(i)}$ for  $\tau \in \TT, i \in \tI$. 
It is known that $\widehat{\WW} 
\cong \TT \widetilde{\times} \widetilde{\WW}$ 
and $\widetilde{\WW} \cong \WW \widetilde{\times} T_{Q^{\vee}}$. 
In particular, the latter isomorphism is given by 
$s_i \mapsto (s_i, id_{\h^{*}})$ for $i \in I$ and 
$s_0 \mapsto (s_{\t}, t_{\t^{\vee}})$, 
where $\t^{\vee}:=\sum_{i \in I} \al^{\vee}_i$. 
The length of an element $\tau w \in \wh{\WW} 
(\tau \in \TT, w \in \widetilde{\WW})$ 
is given by $l_{\wh{\WW}}(\tau w):=l_{\widetilde{\WW}}(w)$, 
where $l_{\widetilde{\WW}}$ is 
the length function of $\widetilde{\WW}$. 

Let $q$ be an indeterminate. 
For $r \in \bbZ$, $m \in \bbN$, 
we define $q$-integers and Gaussian binomial coefficients 
in the rational function field $\bbC(q)$ by 
\begin{eqnarray*}
&&[r]:=\frac{q^{r}-q^{-r}}{q-q^{-1}}, 
\q [m]!:=[m][m-1] \cdots [1], 
\q \left[
\begin{array}{c}
r\\
m
\end{array}
\right]
:=\frac{[r][r-1] \cdots [r-m+1]}{[1][2] \cdots [m]}.
\end{eqnarray*}
Similarly, for $c \in \bbC \, (c \neq 0, \pm 1)$,  we define 
\begin{eqnarray*}
&&[r]_c:=\frac{c^{r}-c^{-r}}{c-c^{-1}}, 
\q [m]_c!:=[m]_c[m-1]_c \cdots [1]_c, 
\q \left[
\begin{array}{c}
r\\
m
\end{array}
\right]_c
:=\frac{[r]_c[r-1]_c \cdots [r-m+1]_c}{[1]_c[2]_c \cdots [m]_c}.
\end{eqnarray*}
We set $[0]!:=[0]_c!:=1$.
\subsection{Definitions}
\newtheorem{def QA}{Definition}[section]
\begin{def QA}
\label{def QA}
The quantum loop algebra $\tU_q:=U_q(\widetilde{\sl}_{n+1})$ 
(resp. the quantum algebra $U_q:=U_q(\sl_{n+1})$, 
the extended quantum algebra $U_q^{'}:=U_q^{'}(\sl_{n+1})$) 
is an associative $\bbC(q)$-algebra generated by 
$\{E_{i}, F_{i}, K_{\mu} \, | \, i \in \tI 
\rm{ (resp. $i \in I$, $i \in I$)},  
\mu \in  \wt{Q} 
\rm{ (resp. $\mu \in Q$, $\mu \in P$)} \}$  
with the relations    
\begin{eqnarray*}
&& K_{\mu} K_{\nu} =K_{\mu+\nu}, \q K_0=1, 
\q K_{\al_0}=K_{\t}^{-1},\\ 
&& K_{\mu} E_{j} K_{\mu}^{-1}=q^{(\mu, \al_j)}E_{j},  
\q K_{\mu} F_{j} K_{\mu}^{-1}=q^{-(\mu, \al_j)}F_{j}, \\ 
&& E_{i} F_{j}- F_{j} E_{i} = \delta_{i,j} 
\frac{K_{\al_i}-K_{\al_i}^{-1}}{q-q^{-1}}, \\
&& \sum_{r=0}^{1-\fra_{ij}} (-1)^r E_{i}^{(r)} E_{j}
E_{i}^{(1-\fra_{i,j}-k)}=
\sum_{r=0}^{1-\fra_{ij}} (-1)^r F_{i}^{(r)} F_{j} 
F_{i}^{(1-\fra_{i,j}-k)}=0
 \q i \neq j, 
\end{eqnarray*}
where 
\begin{eqnarray*}
\q E_{i}^{(r)}:= \displaystyle \frac{1}{[r]!} E_{i}^r, 
\q F_{i}^{(r)} := \displaystyle \frac{1}{[r]!} F_{i}^r 
\q (r \in \bbZ_{+}:=\{0,1,2, \cdots\}). 
\end{eqnarray*}
\end{def QA}
Let $\tU_q^{+}$  (resp. $\tU_q^{-}, \tU_q^{0}$ ) 
be the $\bbC(q)$-subalgebra of $\tU_q$ 
generated by $\{E_{i}\}_{i \in \tI}$ 
(resp.$\{F_{i}\}_{i \in \tI}$, $\{K_{\mu}\}_{\mu \in Q}$). 
Similarly, let $U_q^{+}$  (resp. $U_q^{-}, U_q^{0}$ ) 
be the $\bbC(q)$-subalgebra of $U_q$ 
generated by $\{E_{i}\}_{i \in I}$ 
(resp.$\{F_{i}\}_{i \in I}$, $\{K_{\mu}\}_{\mu \in Q}$). \\
\q It is well known that 
$\tU_q$ (resp. $U_q$) have a Hopf algebra structure 
and its comultiplication is given by 
\begin{eqnarray*}
\D(E_i)=E_i \otimes 1 + K_{\al_i} \otimes E_i, 
 \q \D(F_i)=F_i \otimes K_{\al_i}^{-1} +1 \otimes F_i, 
\q \D(K_{\mu})=K_{\mu} \otimes K_{\mu},
\end{eqnarray*} 
where $i \in \tI$ (resp. $I$), $\mu \in Q$. \\
\q We have a $\bbC$-algebra anti-automorphism 
$\Omega: \tU_q \arr \tU_q$ and 
a $\bbC(q)$-algebra anti-automorphism $\Phi: \tU_q \arr \tU_q$ 
such that 
\begin{eqnarray}
&&\Omega(q)=q^{-1}, \q \Omega(E_i)=F_i, \q \Omega(F_i)=E_i, 
\q \Omega(K_{\mu})=K_{\mu}^{-1},
\label{pro omega} \\
&&\Phi(E_i)=E_i, \q  \Phi(F_i)=F_i, \q  \Phi(K_{\mu})=K_{\mu}^{-1}, 
\label{pro phi}  
\end{eqnarray}
for $i \in \tI, \mu \in Q$.
Let $T_i$ be the $\bbC(q)$-algebra automorphism of $\tU_q$ 
introduced by Lusztig (\cite[\S 37]{L93}: 
\begin{eqnarray}
&& T_i(E_i^{(m)})=(-1)^m q^{-m(m-1)}F_i^{(m)}K_{\al_i}^m, 
\q T_i(F_i^{(m)})=(-1)^m q^{m(m-1)}K_{\al_i}^{-1}E_i^{(m)},  \no\\
&& T_i(E_j^{(m)})=\sum_{r=0}^{-m \fra_{i,j}}(-1)^{r-m\fra_{i,j}}
q^{-r}E_i^{(-m\fra_{i,j}-r)}E_j^{(m)}E_i^{(r)} \q (i \neq j), \no \\
&& T_i(F_j^{(m)})=\sum_{r=0}^{-m \fra_{i,j}}(-1)^{r-m\fra_{i,j}}
q^{r}F_i^{(r)}F_j^{(m)}F_i^{(-m\fra_{i,j}-r)} \q (i \neq j), \no\\
&& T_i(K_{\mu})=K_{s_i(\mu)}, 
\label{pro RA}
\end{eqnarray}
where $i \in \tI, m \in \bbN, \mu \in Q$. 
For $\tau  \in \TT$,
we define $\tU_q$-automorphism $T_{\tau}$ by 
\begin{eqnarray}
T_{\tau}(E_i):=E_{\tau(i)}, \q T_{\tau}(F_i):=F_{\tau(i)}, \q 
T_{\tau}(K_{\al_i}^{\pm 1 }):=K_{\al_{\tau(i)}}^{\pm 1}\,\,
(i \in \tI).
\end{eqnarray}
We obtain that 
\begin{eqnarray}
T_i^{-1}=\Phi T_i \Phi^{-1}, 
\q T_i\Omega=\Omega T_i, 
\q T_{\tau}\Omega=\Omega T_{\tau}. 
\label{pro RA-}
\end{eqnarray} 
Let $w \in \wh{\WW}$ and 
$w=\tau s_{i_1} \cdots s_{i_m} 
(\tau \in \TT, i_1, \cdots , i_m \in \tI, m \in \bbN)$ 
be a reduced expression of $w$. 
Then $T_{w}:=T_{\tau}T_{i_1} \cdots T_{i_m}$ 
is a well-defined $\tU_q$-automorphism, that is, 
$T_{w}$ does not depend on the choice of reduced expression of $w$.
\subsection{Drinfel'd realization}
\q It is known that $\tU_q$ has another realization 
which is called Drinfel'd realization. 
\newtheorem{def DR}[def QA]{Definition}
\begin{def DR}[\cite{D}]
\label{def DR} 
Let $\DD_q$ be an associative $\bbC (q)$-algebra 
generated by 
$\{X_{i,r}^{\pm}, H_{i,s}, K_{\mu} \, | 
\, i \in I, r, s \in \bbZ, s \neq 0, \mu \in Q\}$ 
with the relations 
\begin{eqnarray*}
&&K_{\mu}K_{\nu}=K_{\mu+\nu}, \q  K_0=1, 
\q [K_{\mu}, H_{j,s}]=[H_{i,r}, H_{j,s}]=0, \\
&&K_{\mu}X_{i,r}^{\pm}K_{\mu}^{-1}
=q^{\pm(\mu,\al_i)}X_{i,r}^{\pm}, 
\q [H_{j,s}, X_{i,r}^{\pm}]
=\pm \frac{[r\fra_{i,j}]}{r}X_{i,r+s}, \\
&&X_{i,r+1}^{\pm}X_{j,s}^{\pm}
-q^{\pm \fra_{i,j}}X_{j,s}^{\pm}X_{i,r+1}^{\pm}
=q^{\pm \fra_{i,j}}X_{i,r}^{\pm}X_{j,s+1}^{\pm}
-X_{j,s+1}^{\pm}X_{i,r}^{\pm}, \\
&& [X_{i,r}^{+}, X_{j,s}^{-}]
=\delta_{i,j}\frac{\Psi_{i,r+s}^{+}-\Psi_{i,r+s}^{-}}
{q-q^{-1}}, \\
&&\sum_{\pi \in \SS_{\bf{m}}} \sum_{k=0}^{\bf{m}} (-1)^k 
\left[
\begin{array}{r}
\bf{m}\\
k
\end{array}
\right]
X_{i,r_{\pi(1)}}^{\pm} \cdots X_{i,r_{\pi(k)}}^{\pm} 
X_{j,s}^{\pm} X_{i,r_{\pi(k+1)}}^{\pm} \cdots X_{i,r_{\pi(\bf{m})}}^{\pm}=0, 
\q (i \neq j),
\end{eqnarray*}
for $r_1, \cdots, r_{\bf{m}} \in \bbZ$, 
where $\bf{m}:=1-\fra_{i,j}$, 
$\SS_{\bf{m}}$ is the symmetric group on $\bf{m}$ letters, 
and $\Psi_{i,r}^{\pm}$ are determined by 
\begin{eqnarray*}
\sum_{r=0}^{\infty}\Psi_{i,\pm r}^{\pm}u^{\pm r}
:=K_{\al_i}^{\pm 1}\text{exp}
(\pm (q-q^{-1})\sum_{s=1}^{\infty}H_{i,\pm s}u^{\pm s}),
\end{eqnarray*}
and $\Psi_{i,\pm r}^{\pm}:=0$ if $r<0$.
\end{def DR}
For $i \in I$, let $t_{\L_i^{\vee}}=\tau s_{j_1} \cdots s_{j_r}$ 
($\tau \in \TT$, $j_1, \cdots j_r \in \wt{I}$) 
be a reduced expression of $t_{\L_i^{\vee}}$(see \S 2.1). 
Then we set $T_{\L_i^{\vee}}:=T_{\tau}T_{j_1} \cdots T_{j_r}$.
\newtheorem{thm DR}[def QA]{Theorem}
\begin{thm DR}[\cite{B94a}]
\label{thm DR} 
There exists a $\bbC(q)$-algebra isomorphism 
$T : \DD_q \arr \tU_q$ 
such that 
\begin{eqnarray}
T(X_{i,r}^{+})=(-1)^{ir} T_{\L_i^{\vee}}^{-r}(E_i), 
\q T(X_{i,r}^{-})=(-1)^{ir} T_{\L_i^{\vee}}^{r}(F_i) 
\q (i \in I,r\in\bbZ).\label{sig DRauto}
\end{eqnarray}
\end{thm DR}
In particular, by \cite{B94a} \S4 Remark and \cite{CP94a} \S2.5, 
we obtain the inverse map of $T$:
\begin{eqnarray}
&&T^{-1}(E_i)=X_{i,0}^{+}, \q T^{-1}(F_i)=X_{i,0}^{-}, 
\q T^{-1}(K_{\mu})=K_{\mu}, \no \\
&&T^{-1}(E_0)=(-1)^{m+1}q^{n+1}[X_{n,0}^{-},  \cdots [
X_{m+1,0}^{-},[ X_{1,0}^{-},  \cdots [X_{m-1,0}^{-},
X_{m,1}^{-}]_{q^{-1}} \cdots ]_{q^{-1}}K_{\t}^{-1}, \no \\
&&T^{-1}(F_0)=(-1)^{m+n}[X_{n,0}^{+},  \cdots [
X_{m+1,0}^{+},[ X_{1,0}^{+},  \cdots [X_{m-1,0}^{+},
X_{m,-1}^{+}]_{q^{-1}} \cdots ]_{q^{-1}}K_{\t}, 
\label{pro DR}
\end{eqnarray}
for $m, i \in I$, 
where $[u,v]_{q^{\pm1}}:=uv-q^{\pm1}vu$ for $u,v \in \tU_q$ 
($T$ is independent of the choice of $m$). 
We idntify $\DD_q$ with $\tU_q$ by this isomorphism $T$.
\subsection{PBW basis}
\q Let $w_0$ be a longest element in $\WW$ 
and $w_0=s_{i_1} \cdots s_{i_N}$ be a reduced expression of $w_0$. 
We set 
$\ga_1:=\al_{i_1},  \ga_2:=s_{i_1}(\al_{i_2}), 
 \cdots, \ga_N:=s_{i_1} \cdots s_{i_{N-1}}(\al_{i_N})$.
By the theory of the classical Lie algebra, 
we have $\D_{+}=\{\ga_1, \cdots, \ga_N\}$. 
Define the root vectors in $U_q$ by 
\begin{eqnarray} 
\bar{E}_{\ga_k}:=T_{i_1} \cdots T_{i_{k-1}}(E_{\al_{i_k}}), 
\q \bar{F}_{\ga_k}:=\Om (E_{\ga_k}),
\label{def root vectors}
\end{eqnarray}
for $1 \leq k \leq N$, 
where $E_{\al_i}:=E_i$, $F_{\al_i}:=F_i$ $(i \in I$). 
We set 
\begin{eqnarray} 
&&\bbZ_{+}^{\D_{+}}:=\{c: \D_{+} \arr \bbZ_{+}; \textrm{map}\}, 
\q B_q^{0}:=\{K_{\mu} \, | \, \mu \in Q\},  \no\\ 
&&B_q^{+}:=\{\prod _{\ga \in \D_{+}}^{<} \bar{E}_{\ga}^{c(\ga)} 
\, | \, c \in \bbZ_{+}^{\D_{+}}\}, 
\q B_q^{-}:=\Om(B_q^{+}), 
\q B_q:=B_q^{-}B_q^{0}B_q^{+}, 
\label{def bq}
\end{eqnarray}
where $<$ means that the product is ordered by 
$\ga_1 < \cdots <\ga_N$. 
\newtheorem{thm PBW}[def QA]{Theorem}
\begin{thm PBW}[\cite{L90a} \S1]
\label{thm PBW}
$B_q^{\star}$ (resp. $B_q$) is a $\bbC(q)$-basis of $U_q^{\star}$ 
(resp. $U_q$) for $\star \in \{-, 0, +\}$.
\end{thm PBW}
The following notations and facts are given in \cite{BCP}. 
Let $\rho^{\vee}:=\frac{1}{2} \sum_{\al \in \D_+} \al^{\vee} 
( \in P^{\vee})$ 
and $t_{2\rho^{\vee}}=s_{j_1} \cdots s_{j_{\widetilde{N}}}$ 
be a reduced expression of $t_{2\rho^{\vee}} (2 \rho^{\vee} \in
Q^{\vee})$. 
Define the doubly infinite sequence $(\cdots, i_{-1}, i_0, i_1, \cdots)$ 
by $i_k:=j_{k^{'}}$ if $k \equiv k^{'}$ (mod $\widetilde{N}$) 
for all $k \in \bbZ, 1 \leq k^{'} \leq \widetilde{N}$. 
We set 
\begin{eqnarray*} 
\beta_k:=s_{i_0}s_{i_{-1}} \cdots s_{i_{k+1}}(\al_{i_k}) \q (k \leq 0), 
\q \beta_k:=s_{i_1}s_{i_2} \cdots s_{i_{k-1}}(\al_{i_k}) \q (k > 0). 
\end{eqnarray*} 
Then we have $\Dr_{+}=\{\b_k\}_{k \in \bbZ}$. 
Define a total order on $\widetilde{\D}_{+}(I)$ by  
\begin{eqnarray}
\hspace{-30pt}\b_0 < \b_{-1}< \b_{-2}< \cdots <(1, \d)< \cdots <(n,\d) 
<(1, 2\d)< \cdots< (n,2\d)< \cdots < \b_2 < \b_1. 
\label{sig TO}
\end{eqnarray}
We set $E_{\al_i}:=E_i$, $F_{\al_i}:=F_i$ $(i \in \wt{I}$). 
 Define the positive real root vectors in $\tU_q$ by 
\begin{eqnarray}
E_{\b_k}:=T_{i_0}^{-1}T_{i_{-1}}^{-1} \cdots T_{i_{k+1}}^{-1}
(E_{\al_{i_k}}) \q (k \leq 0), 
\q E_{\b_k}:=T_{i_1}T_{i_{2}} \cdots T_{i_{k-1}}
(E_{\al_{i_k}}) \q (k > 0), 
\label{def ARV} 
\end{eqnarray}
and the positive imaginary root vectors $E_{(i,r \d)}$ by 
\begin{eqnarray}
\rm{exp}( (q-q^{-1}) \sum_{s=1}^{\infty} E_{(i,s\d)}u^k)
:=1 +\sum_{s=1}^{\infty}(q-q^{-1})\wh{E}_{(i,s\d)}u^k, 
\label{def IRV}
\end{eqnarray}
where $\wh{E}_{(i,s\d)}:=E_{-\al_i+s\d}E_{i}-q^{-2}E_iE_{-\al_i+s\d}$ 
for $i \in I$, $s \in \bbN$. 
Define the negative root vectors by 
$F_{\b}:=\Om(E_{\b})$ for $\b \in \widetilde{\D}_{+}(I)$. 
We set 
\begin{eqnarray}
&&\bbZ^{\tD_{+}(I)}_{+}:=\{c: \tD_{+}(I) \arr \bbZ_{+} ;
\textrm{map} \, |\, \#\{c(\b) \neq 0\} < \infty\}, 
\q  \tB_q^{0}:=\{K_{\mu} \, |\, \mu \in Q\},  \no \\ 
 &&\tB_q^{+}:=\{\prod_{\b \in \tD_{+}(I)}^{<} E_{\b}^{c(\b)} \, 
| \, c \in \bbZ^{\tD_{+}(I)}_{+}\}, 
\q  \tB_q^{-}:=\Om(\tB_q^{+}), 
\q \tB_q:=\tB_q^{-}\tB_q^{0}\tB_q^{+}, 
\label{sig APBW}
\end{eqnarray}
where $<$ is the total order as in (\ref{sig TO}). 
\newtheorem{thm APBW}[def QA]{Theorem}
\begin{thm APBW}[\cite{B94b}]
\label{thm APBW}
$\tB_q^{\star}$ (resp. $\tB_q$) is a $\bbC(q)$-basis of $\tU_q^{\star}$ 
(resp. $\tU_q$) for $\star \in \{-, 0, +\}$.
\end{thm APBW}
By \cite{BCP} Lemma 1.5, we obtain 
\begin{eqnarray}
&&X_{i,r}^{+}=(-1)^{ir}E_{\al_i+r\d} \q (r \geq 0), 
\q X_{i,r}^{+}=(-1)^{ir-1}F_{-\al_i-r\d}K_i^{-1} \q (r <0), \no \\
&&X_{i,r}^{-}=(-1)^{ir-1}K_iE_{-\al_i+r\d} \q (r > 0), 
\q X_{i,r}^{-}=(-1)^{ir}F_{\al_i-r\d} \q (r  \leq0), \no \\
&&H_{i,s}=(-1)^{is}E_{(i,s\d)}, 
\q \Psi_{i,s}^{+}=(-1)^{is}(q-q^{-1})K_i \wh{E}_{(i,s\d)}, 
\label{rem DRV}
\end{eqnarray}
for $i \in I$, $r \in \bbZ, s \in \bbN$.
\subsection{Evaluation homomorphisms}

There exists a $\bbC(q)$-algebra  homomorphism
 $U_q(\sl_{n+1}) \arr U_q(\widetilde{\sl}_{n+1})$ such that 
\begin{eqnarray}
E_i \mapsto X_{i,0}^{+}, \q F_i \mapsto X_{i,0}^{-}, 
\q K_{\mu} \mapsto K_{\mu}, 
\label{pro embed}
\end{eqnarray}
for $i \in I, \mu \in Q$. 
Moreover, for $m \in I$, $0 \leq k \leq n-m$, 
there exists a $\bbC(q)$-algebra  homomorphism
$U_q(\widetilde{\sl}_{m+1}) \arr 
U_q(\widetilde{\sl}_{n+1})$ such that 
\begin{eqnarray}
X_{i,r}^{\pm} \mapsto X_{i+k,r}^{\pm}, \q 
H_{i+k,s} \mapsto H_{i,s}, 
\q K_{\al_i} \mapsto K_{\al_{i+k}}, 
\label{pro embed2}
\end{eqnarray}
for $1 \leq i \leq m$, $r,s \in \bbZ(s \neq 0)$. 
Hence we can regard any $U_q(\widetilde{\sl}_{n+1})$-module 
as a $U_q(\sl_{n+1})$-module 
and $U_q(\widetilde{\sl}_{m+1})$-module. 
Let $U^{'}_q(\sl_{n+1})$ be the extended quantum algebra 
in Definition \ref{def QA}.
By the following proposition, we can regard any 
$U_q^{'}(\sl_{n+1})$-module as a 
$U_q(\widetilde{\sl}_{n+1})$-module. 
\newtheorem{pro EH}[def QA]{Proposition}
\begin{pro EH}[\cite{J} \S2, \cite{CP94a} Proposition 3.4]
\label{pro EH}
For any $\bf{a} \in \bbC^{\times}$, 
there exist $\bbC(q)$-algebra homomorphisms 
$ev_{\bf{a}}^{\pm}:U_q(\widetilde{\sl}_{n+1})\arr U^{'}_q(\sl_{n+1})$ 
such that 
\begin{eqnarray*}
&&ev_{\bf{a}}^{\pm}(E_i)=E_i, 
\q ev_{\bf{a}}^{\pm}(F_i)=F_i,
\q ev_{\bf{a}}^{\pm}(K_{\mu})=K_{\mu}, \\
&&ev_{\bf{a}}^{+}(E_0)
=q^{-1}\bf{a}K_{\L_1}K_{\L_n}^{-1}[F_n, [F_{n-1}, \cdots, [
F_2,F_1]_{q^{-1}}  \cdots ]_{q^{-1}}, \\
&&ev_{\bf{a}}^{-}(E_0)
=q^{-1}\bf{a}K_{\L_1}^{-1}K_{\L_n}[F_1, [F_{2}, \cdots, [
F_{n-1},F_n]_{q^{-1}}  \cdots ]_{q^{-1}}, \\ 
&&ev_{\bf{a}}^{+}(F_0)=(-1)^{n-1}q^{n}{\bf{a}}^{-1}K_{\L_1}^{-1}K_{\L_n}[E_n, [
E_{n-1}, \cdots, [E_2,E_1]_{q^{-1}}  \cdots ]_{q^{-1}}, \\
&&ev_{\bf{a}}^{-}(F_0)=(-1)^{n-1}q^{n}{\bf{a}}^{-1}K_{\L_1}K_{\L_n}^{-1}[E_1, [
E_{2}, \cdots, [E_{n-1},E_n]_{q^{-1}}  \cdots ]_{q^{-1}}, 
\end{eqnarray*}
\end{pro EH}
for $i \in I$ and $\mu\in P$
By (\ref{pro RA}), (\ref{pro RA-}), we obtain 
\begin{eqnarray}
&&ev_{\bf{a}}^{+}(E_0)=q^{-n}\bf{a}
K_{\L_1}K_{\L_n}^{-1}T_1^{-1} \cdots T_{n-1}^{-1}(F_n), 
\q ev_{\bf{a}}^{-}(E_0)=q^{-n}\bf{a}
K_{\L_1}^{-1}K_{\L_n}T_n^{-1} \cdots T_{2}^{-1}(F_1), 
\no  \\
&&ev_{\bf{a}}^{+}(F_0)=q^{n}{\bf{a}}^{-1}K_{\L_1}^{-1}K_{\L_n}T_1^{-1} 
\cdots T_{n-1}^{-1}(E_n),
\q ev_{\bf{a}}^{-}(F_0)=q^{n}{\bf{a}}^{-1}K_{\L_1}K_{\L_n}^{-1}T_n^{-1} 
\cdots T_{2}^{-1}(E_1). \no \\
\label{pro EH2}
\end{eqnarray}

\section{Quantum algebras at roots of unity}
\setcounter{equation}{0}
\renewcommand{\theequation}{\thesection.\arabic{equation}}
\q In the rest of this paper, we fix the following notations. 
Let $l$ be an odd integer greater than $2$ and 
$\e$ a primitive $l$-th root of unity. 
Moreover, we assume g.c.d$(l,n+1)=1$. 
By \cite{BK} Lemma 2.1 and Corollary 2.1, 
we obtain that g.c.d$(l,n+1)=1$ if and only if 
$\rm{det}([k\fra_{i,j}])_{i,j \in I} \neq 0$ 
for any $k \in \bbZ$ such that $k \not\equiv 0$ mod $l$. 

\subsection{Quantum algebras of non-restricted type} 
\q Let $\AA:=\bbC[q,q^{-1}]$ be the Laurent polynomial ring 
and $\tu_{\AA}$ (resp. $U_{\AA}$) 
the $\AA$-subalgebra of $\tu_q$ (resp. $U_q$) generated by 
$\{E_i, F_i, K_{\al_j}, [K_{\al_j};0] \, | \, 
i \in \tI (\rm{resp. } i \in I), j \in I \}$, 
where $[K_{\al_j} ; 0]:=(K_{\al_j}-K_{\al_j}^{-1})/(q-q^{-1})$ 
for $j \in I$. 
Let $\tU_{\AA}^{+}$ (resp. $\tU_{\AA}^{-}, \tU_{\AA}^{0}$) 
be the $\AA$-subalgebra of $\tU_{\AA}$ generated by 
$\{E_i\}_{i \in \tI}$ (resp. $\{F_i\}_{i \in \tI}$, 
$\{K_{\al_i}, [K_{\al_i};0] \}_{i \in I}$) 
and $U_{\AA}^{+}$ (resp. $U_{\AA}^{-}, U_{\AA}^{0}$) 
the $\AA$-subalgebra of $U_{\AA}$ generated by 
$\{E_i\}_{i \in I}$ (resp. $\{F_i\}_{i \in I}$, 
$\{K_{\al_i}, [K_{\al_i};0] \}_{i \in I}$). 
We have triangular decompositions 
$\tU_{\AA}=\tU_{\AA}^{-}\tU_{\AA}^{0}\tU_{\AA}^{+}$ and 
$U_{\AA}=U_{\AA}^{-}U_{\AA}^{0}U_{\AA}^{+}$. We set 
\begin{eqnarray*}
\tB_{\AA}^{+}:=\tB_{q}^{+}, \, \, \tB_{\AA}^{-}:=\tB_{q}^{-}, 
\, \,  \tB_{\AA}^{0}:=\{\prod_{i \in I}K_{\al_i}^{\d_i}[K_{\al_i};0]^{m_i} \, |\, 
m_i \in \bbZ_{+}, \d_i=0 \rm{ or }1\}, 
\, \, \tB_{\AA}:=\tB_{\AA}^{-}\tB_{\AA}^{0}\tB_{\AA}^{+}, \\
B_{\AA}^{+}:=B_{q}^{+}, \, \, B_{\AA}^{-}=B_{q}^{-}, 
\, \, B_{\AA}^{0}:=\{\prod_{i \in I}K_{\al_i}^{\d_i}[K_{\al_i};0]^{m_i} \, |\, 
m_i \in \bbZ_{+}, \d_i=0 \rm{ or }1\}, 
\, \, B_{\AA}:=B_{\AA}^{-}B_{\AA}^{0}B_{\AA}^{+}.
\end{eqnarray*}
\q We have $T_i(\tU_{\AA}) \subset \tU_{\AA}$ by (\ref{pro RA}) 
and $T_i^{-1}(\tU_{\AA}) \subset \tU_{\AA}$ by (\ref{pro RA-}) 
for $i \in I$. 
Hence, by (\ref{def ARV}) and (\ref{rem DRV}), 
we have $E_{\b}, F_{\b}, X_{i,r}^{\pm}, H_{i,s} \in \tU_{\AA}$ 
for all $\b \in \tD_{+}(I)$, $i \in I$, $r,s \in \bbZ(s \neq 0)$. 
Similarly, we obtain 
$T_i^{\pm 1}(U_{\AA}) \subset U_{\AA}$ and 
$\bar{E}_{\ga}, \bar{F}_{\ga} \in U_{\AA}$ 
for all $\ga \in \D_{+}$ by (\ref{def root vectors}). 
Thus we obtain 
$\tB_{\AA}^{\star}, \tB_{\AA} \subset \tU_{\AA}$ and 
$B_{\AA}^{\star}, B_{\AA} \subset U_{\AA}$ ($\s=$).
\newtheorem{sec 3}{Proposition}[section]
\begin{sec 3} 
$B_{\AA}^{\star}$ (resp. $B_{\AA}$) is a $\AA$-basis of 
$U_{\AA}^{\star}$ (resp. $U_{\AA}$) for $\s=$.
\label{sec 3}
\end{sec 3}
Proof. 
By Theorem \ref{thm PBW}, $B_{q}^{\star}$ (resp. $B_{q}$) 
is $\bbC$-linearly independent in $U_{q}^{\star}$ 
(resp. $U_{q}$) for $\s=$. 
Thus $B_{\AA}^{\star}$ (resp. $B_{\AA}$) 
is $\AA$-linearly independent in $U_{\AA}^{\star}$ (resp. $U_{\AA}$). 
Let $V_{\AA}^{\star}$ (resp. $V_{\AA}$) be the $\AA$-subalgebra of
$U_{\AA}^{\star}$ (resp. $U_{\AA}$) generated by 
$B_{\AA}^{\star}$ (resp. $B_{\AA}$). 
It is enough to prove that 
$U_{\AA}^{\star}V_{\AA}^{\star} \subset V_{\AA}^{\star}$ for all $\s=$. 
Indeed, if we can prove this claim, then we obtain 
$U_{\AA}^{\star}=V_{\AA}^{\star}$
and $U_{\AA}=U_{\AA}^{-}U_{\AA}^{0}U_{\AA}^{+}
=V_{\AA}^{-}V_{\AA}^{0}V_{\AA}^{+}=V_{\AA}$. 
So $U_{\AA}^{\star}$ (resp. $U_{\AA}$) is generated by $B_{\AA}^{\star}$ 
(resp. $B_{\AA}$) as $\AA$-module. \\
\q By the following formula, 
we have $K_{\al_i}^{\pm 1}(\prod_{i \in I}K_{\al_i}^{\d_i}[K_{\al_i};0]^{m_i}) 
\in V_{\AA}^{0}$ for all $i \in I$, $m_i \in \bbZ_{+}$:
\begin{eqnarray*}
&& K_{\al_i}^2=K_{\al_i}(K_{\al_i}-K_{\al_i}^{-1})+1=(q-q^{-1})K_{\al_i}[K_{\al_i};0]+1 \in V_{\AA}^{0}, \\
&& K_{\al_i}^{-1}=K_{\al_i}-(K_{\al_i}-K_{\al_i}^{-1})=K_{\al_i}-(q-q^{-1})[K_{\al_i};0] \in V_{\AA}^{0}. 
\end{eqnarray*}
Thus we obtain $U_{\AA}^{0}V_{\AA}^{0} \subset V_{\AA}^{0}$. 
By \cite{DK} Lemma 1.7, we get the following formula: 
for $\al, \b \in \D_{+}$ such that $\b>\al$,  
\begin{eqnarray*}
 E_{\b}E_{\al}=\sum_{c \in \bbZ_{+}^{\D_{+}}}a_c
\prod_{\ga \in \D_{+}}^{<}E_{\ga}^{c(\ga)},
\end{eqnarray*}
where $a_c \in \AA$. So we obtain the case of $\star=+$. 
Similarly, by using the automorphism $\Om$ (see (\ref{pro omega})), 
we obtain the case of $\star=-$. 
\qed \\ \\
\q By \cite{BK} Proposition 1.7(c), we obtain the following formula: 
for $\al, \b \in \tD_{+}(I)$ such that $\b >\al$, 
\begin{eqnarray*}
 E_{\b}E_{\al}=q^{(\al,\b)}E_{\al}E_{\b}
+\sum_{\al <\ga_1 < \cdots <\ga_m <\b}c_{\ga}E_{\ga_1}^{a_1}
\cdots E_{\ga_m}^{a_m},
\end{eqnarray*}
where $c_{\ga} \in \bbC[q,q^{-1}]$ 
for $\ga=(\ga_1, \cdots , \ga_m) \in \tD_{+}(I)^m$. 
So, by the similar way to the proof of Proposition \ref{sec 3}, 
we obtain the following proposition.
\newtheorem{pro AAPBW}[sec 3]{Proposition}
\begin{pro AAPBW} 
\label{pro AAPBW}
$\tB_{\AA}^{\star}$ (resp. $\tB_{\AA}$) is a $\AA$-basis of 
$\tU_{\AA}^{\star}$ (resp. $\tU_{\AA}$) for $\s=$.
\end{pro AAPBW}
Now we define the quantum algebras of non-restricted type.
We regard $\bbC$ as $\AA$-module by 
$g(q).c:=g(\e)c$ for $g(q) \in \AA$, $c \in  \bbC$ 
and denote it by $\bbC_{\e}$. 
We define 
\begin{eqnarray*}
 \tu_{\e}:=\tu_{\AA} \otimes_{\AA} \bbC_{\e} 
\q (\rm{resp. }U_{\e}:=U_{\AA} \otimes_{\AA} \bbC_{\e}).
\end{eqnarray*}
Then we call $\tu_{\e}$ (resp. $U_{\e}$) 
``quantum loop algebra (resp. quantum algebra) of non-restricted
type (or De Concini-Kac type)'' (see \cite{BK}, \cite{DK}). 
For $\s=$, we set 
$\tU_{\e}^{\star}:=\tU_{\AA}^{\star} \otimes_{\AA} 1$ 
(resp. $U_{\e}^{\star}:=U_{\AA}^{\star} \otimes_{\AA} 1$). 
We simply denote $u \otimes 1$ 
by $u$ for $u \in \tu_{\AA}$ (resp. $U_{\AA}$). 
\newtheorem{rem NS}[sec 3]{Remark}
\begin{rem NS} 
\label{rem NS}
(a) In \cite{BK} (resp. \cite{DK}), 
$\tu_{\e}$ (resp. $U_{\e}$) is defined by 
$\tu_{\e}:=\tu_{\AA^{\e}} / (q-\e)\tu_{\AA^{\e}}$ 
(resp. $U_{\e}:=U_{\AA} / (q-\e)U_{\AA}$), 
where $\AA^{\e}:=\{g(q) \in \bbC(q) \, | \, 
g(q) \rm{ has no poles at } q=\e\} (\supset \AA)$ 
and $(q-\e)\tu_{\AA^{\e}}$ (resp. $(q-\e)U_{\AA}$) 
is the two-sided ideal 
of $\bbC$-algebra $\tu_{\AA^{\e}}$ 
(resp. $U_{\AA}$) generated by $(q-\e)$. 
But, by the universality of tensor product, we obtain
\begin{eqnarray*}
&&\tU_{\AA} \otimes_{\AA}\bbC_{\e} \cong \tU_{\AA}/(q-\e)\tU_{\AA} \cong 
\tU_{\AA^{\e}}/(q-\e)\tU_{\AA^{\e}} 
\rm{  (as $\bbC$-algebra)}, \no \\
&& U_{\AA} \otimes_{\AA}\bbC_{\e} \cong U_{\AA}/(q-\e)U_{\AA} 
 \cong U_{\AA^{\e}}/(q-\e)U_{\AA^{\e}}
\rm{  (as $\bbC$-algebra)}.
\end{eqnarray*}
\q (b) $\tu_{\e}$ (resp. $U_{\e}$) is the associative algebra over $\bbC$ 
on generators $\{E_i, F_i, K_{\mu}$ $|$ $i \in \tI \, 
(\rm{resp. } I), \, \mu \in \wt{Q} \, (\rm{resp. }Q)\}$
and defining relations in Definition \ref{def QA} 
replaced $q$ by $\e$ 
(see \cite{BK} \S 1.9 and \cite{DK} \S 1.5).
\end{rem NS}
We set $\tB_{\e}^{\star}:=\tB_{\AA} \otimes _{\AA} 1$ $(\s=)$. 
Similarly, we define $\tB_{\e}, B_{\e}^{\star}$ and $B_{\e}$. 
\newtheorem{lem basis}[sec 3]{Lemma}
\begin{lem basis} 
\label{lem basis}
Let $\{v_j\}_j$ be a $\AA$-basis in $\tU_{\AA}$ (resp. $U_{\AA}$). 
Then $\{v_j+(q-\e)\tU_{\AA}\}_j$ (resp. $\{v_j+(q-\e)U_{\AA}\}_j$ )
is a $\bbC$-basis of $\tU_{\AA}/(q-\e)\tU_{\AA}$ 
(resp. $U_{\AA}/(q-\e)U_{\AA}$).
\end{lem basis}
Proof. 
$\tU_{\AA}/(q-\e)\tU_{\AA}$ is spanned by 
$\{v_j+(q-\e)\tU_{\AA}\}_j$ as $\bbC$-vector space. 
So we shall prove that 
$\{v_j+(q-\e)\tU_{\AA}\}_j$ is linearly independent over $\bbC$ 
in $\tU_{\AA}/(q-\e)\tU_{\AA}$.
We assume that 
$\sum_jc_j(v_j+(q-\e)\tU_{\AA})=0$ 
($c_j \in \bbC$, $\#\{j \, |\, c_j \neq 0\}<\infty$). 
Then $\sum_{j}c_jv_j \in (q-\e)\tU_{\AA}$. 
Since $\tU_{\AA}$ is generated by $\{v_j\}_j$ as $\AA$-algebra, 
there exist $c_{j,m} \in \bbC$ 
($m \in \bbZ$, $\#\{(j,m) \, |\, c_{j,m} \neq 0\}<\infty$) such that 
$\sum_jc_jv_j=(q-\e)\sum_{j,m}c_{j,m}q^m v_j$ in $\tU_{\AA}$. 
Since $\{v_j\}_j$ is linearly independent over $\AA$ in $\tU_{\AA}$, 
we obtain $c_j=(q-\e)\sum_{m \in \bbZ}c_{j,m}q^m$ for all $j$. 
Therefore we obtain $c_j=c_{j,m}=0$ for any $j$ and $m$. 
Similarly, we obtain the case of $U_{\AA}$.
\qed \\ \\
\q By this lemma, Proposition \ref{pro AAPBW}
and Remark \ref{rem NS}(a), we obtain the following proposition.
\newtheorem{pro ANRPBW}[sec 3]{Proposition}
\begin{pro ANRPBW} 
\label{pro ANRPBW} 
$\tB_{\e}^{\star}$ (resp. $\tB_{\e}$) 
is a $\bbC$-basis of $\tU_{\e}^{\star}$ (resp. $\tU_{\e}$) 
for $\star \in \{-, 0, +\}$.
\end{pro ANRPBW}
The classical case of this proposition is given in \cite{DK} \S 1.7.
\newtheorem{pro NRPBW}[sec 3]{Proposition}
\begin{pro NRPBW} 
\label{pro NRPBW}
$B_{\e}^{\star}$ (resp. $B_{\e}$) 
is a $\bbC$-basis of $U_{\e}^{\star}$ (resp. $\ue$)
for $\star \in \{-, 0, +\}$.
\end{pro NRPBW}
Let $Z(\tU_{\e})$ (resp. $Z(\ue)$) be the center of 
$\tue$ (resp. $\ue$) 
and $\wt{Z}_0$ (resp. $Z_0$) be the $\bbC$-subalgebra of $\tUe$ (resp. $\ue$) 
generated by 
$\{E_{\b}^l, F_{\b}^l, E_{(i,sl\d)}, F_{(i,sl\d)}, K_{\mu}^l 
\, | \, \b \in \wt{\D}_{+}^{\rm{re}}, 
i \in I, s \in \bbN,  \mu \in Q\}$ 
(resp. $\{\bar{E}_{\ga}^l, \bar{F}_{\ga}^l, K_{\mu}^{l} \, | \, 
\ga \in \D_{+}, \mu \in Q\}$).
\newtheorem{pro ACE}[sec 3]{Proposition}
\begin{pro ACE}[\cite{BK} Lemma 2.2, Proposition 2.3]
\label{pro ACE}
We have $\wt{Z}_0=Z(\tU)$.
\end{pro ACE}
\newtheorem{pro CE}[sec 3]{Proposition}
\begin{pro CE}[(\cite{DK} Corollary 3.1]
\label{pro CE}
We have $Z_0 \subset Z(\ue)$. 
\end{pro CE}
For $m \in \bbN$, 
we set $\bbZ_m:=\{0,1, \cdots, m-1\} \subset \bbZ_{+}$ and 
$Q_m:=\bigoplus_{i \in I} \bbZ_m\al_i$. 
Let $I_{\e}$ be the two sided ideal of $\ue$ 
generated by $\{\bar{E}_{\ga}^l, \bar{F}_{\ga}^l, K_{\mu}^{2l}-1 
\, | \, \ga \in \D_{+}, \mu \in Q\}$ 
and set $(\ue / I_{\e})^{\star}:=\{u+I_{\e} \, | \, u \in \ue^{\star}\}
\subset \ue/I_{\e}$ for $\s=$. 
We set 
\begin{eqnarray}
&&\bbZ_l^{\D_{+}}:=\{c:\D_{+} \arr \bbZ_l ; \rm{map}\},  
 \q B_l^{+}:=\{(\prod_{\ga \in \D_{+}}^{<}\bar{E}_{\ga})+I_{\e} \, 
| \, \ga \in \bbZ_l^{\D_{+}}\}, \no \\
&& B_l^{-}:=\Om(B_l^{+}),  
\q B_l^{0}:=\{K_{\mu}+I_{\e} \, | \, \mu \in Q_{2l}\}, 
\q B_l:=B_l^{-}B_l^{0}B_l^{+}.
\label{def blpbw}
\end{eqnarray}
\newtheorem{pro lPBW}[sec 3]{Proposition}
\begin{pro lPBW}
\label{pro lPBW}
$B_l^{\star}$ (resp. $B_l$) is a $\bbC$-basis of 
$(\ue/I_{\e})^{\star}$ (resp. $\ue/I_{\e}$) for $\s=$. 
\end{pro lPBW}
Proof. 
We shall prove that $B_l^{+}$ is a $\bbC$-basis of $(\ue/I_{\e})^{+}$. 
We can also prove the other cases similarly. 
Let $V_l^{+}$ be the $\bbC$-subspace of $(\ue/I_{\e})^{+}$ 
spanned by $B_{l}^{+}$ 
and $u+I_{\e} \in (\ue/I_{\e})^{+}$. 
By Theorem \ref{thm PBW}, 
there exist $a_c \in \bbC (c \in \bbZ_{+}^{\D_{+}})$ 
such that $u+I_{\e}=\sum_{c \in \bbZ_{+}^{\D_{+}}}a_c 
\prod_{\ga \in \D_{+}}^{<} \bar{E}_{\ga}^{c(\ga)}+I_{\e}$. 
If $c$ is an element in $\bbZ_{+}^{\D_{+}}\verb|\|  \bbZ_l^{\D_{+}}$, 
then there exists $\ga_0 \in \D_{+}$ such that 
$c(\ga_0) \geq l$. 
Then we have $\prod_{\ga \in \D_{+}}^{<} \bar{E}_{\ga}^{c(\ga)} \in I_{\e}$. 
Hence 
$ u+I_{\e}=\sum_{c \in \bbZ_l^{\D_{+}}}a_c 
\prod_{\ga \in \D_{+}}^{<} \bar{E}_{\ga}^{c(\ga)}+I_{\e} \in V_l^{+}$. 
So $(U_{\e}/I_{\e})^{+}$ is spanned by $B_l^{+}$. \\
\q Let $u=\sum_{c \in \bbZ_l^{\D_{+}}}a_c 
\prod_{\ga \in \D_{+}}^{<} \bar{E}_{\ga}^{c(\ga)} \in U_{\e}^{+}$ 
$(a_c \in \bbC, c \in \bbZ_l^{\D_{+}}$). 
We assume $u+I_{\e}=0$ in $(\ue /I_{\e})^{+}$. 
Then we have $u \in I_{\e} \cap \ue^{+}$.
By Theorem \ref{thm PBW}, we obtain 
\begin{eqnarray*}
 I_{\e} \cap \ue^{+}
=(\sum_{\ga \in \D_{+}}\ue \bar{E}_{\ga}^{l}
+\sum_{\mu \in Q}\ue(K_{\mu}^{2l}-1)
+\sum_{\ga \in \D_{+}}\ue \bar{F}_{\ga}^l) 
\cap \ue^{+}
=\sum_{\ga \in \D_{+}}\ue^{+}\bar{E}_{\ga}^l.
\end{eqnarray*}
Hence there exists 
$u^{'}=\sum_{c^{'} \in \bbZ_{+}^{\D_{+}}}b_{c^{'}} 
\prod_{\ga \in \D_{+}}^{<} \bar{E}_{\b}^{c^{'}(\b)} \in \ue^{+}
(b_{c^{'}} \in \bbC, c^{'} \in \bbZ_l^{\D_{+}})$ 
such that $u=\sum_{\ga \in \D_{+}} u^{'}\bar{E}_{\ga}^{l}$. 
Since $\bar{E}_{\ga}^l$ is a central element in $\ue$ 
(see Proposition \ref{pro CE}), we have 
\begin{eqnarray*}
 \sum_{\ga \in \D_{+}} u^{'}\bar{E}_{\ga}^{l} 
\in \sum_{c \in \bbZ_{+}^{\D_{+}} \verb|\| \bbZ_l^{\D_{+}} }
\bbC (\prod_{\ga \in \D_{+}}^{<} \bar{E}_{\ga}^{c(\ga)}).
\end{eqnarray*} 
Thus, by Proposition \ref{pro NRPBW}, we get
$a_c=0$ for any $c \in \bbZ_l^{\D_{+}}$. 
Therefore $B_l^{+}$ is linearly independent in 
$(U_{\e}/I_{\e})^{+}$.
\qed \\ \\
\q Let $\tI_{\e}$ be the two sided ideal of $\tue$ 
generated by 
$\{E_{\b}^l, F_{\b}^l, E_{(i,sl\d)}, F_{(i,sl\d)}, K_{\mu}^{2l}-1 \,
| \, \b \in \tD_{+}^{\rm{re}}, i \in I, s \in \bbN, \mu \in Q\}$ 
and set 
$(\tue/\tI_{\e})^{\star}:=\{u+\tI_{\e} \, | \, u \in \tue^{\star}\}$ 
for $\s=$. 
We define  
\begin{eqnarray}
&&\bbZ_l^{\tD_{+}(I)}:=\{c \in \bbZ^{\tD_{+}(I)}_{+} \, | \, 
c(\b) \in \bbZ_l, c((i,sl\d))=0 
\, (\b \in \tD_{+}^{\rm{re}}, i \in I, s \in \bbN) \}, \no \\ 
&& \tB_l^{+}:=\{\prod_{\b \in \tD_{+}(I)}^{<}E_{\b}^{c(\b)} + \tI_{\e} \, 
| \, c \in \bbZ_l^{\tD_{+}(I)}\}, 
\q \tB_l^{-}:=\Om(\tB_l^{+}), \no \\
&& \tB_l^0:=\{K_{\mu}+\tI_{\e} \, | \, \mu \in Q_{2l}\}, 
\q \tB_l:=\tB_l^{-}\tB_l^{0}\tB_l^{+}.
\label{def ARSl-basis}
\end{eqnarray} 
Then, by Theorem \ref{thm APBW} and Proposition \ref{pro ACE}, 
we obtain the following proposition. 
\newtheorem{pro lAPBW}[sec 3]{Proposition}
\begin{pro lAPBW}
\label{pro lAPBW}
$\tB_l^{\star}$ (resp. $\tB_l$) is a $\bbC$-basis of 
$(\tue/\tI_{\e})^{\star}$ (resp. $(\tue/\tI_{\e})$) for $\s=$. 
\end{pro lAPBW}
The proof of this proposition is similar to the one of 
Proposition \ref{pro lPBW}. 
\subsection{Quantum algebras of restricted type}
\q Let $\tU_{\AA}^{\rm{res}}$ (resp. $U_{\AA}^{\rm{res}}$) 
be the $\AA$-subalgebra of $\tU_q$ (resp. $U_q$) 
generated by 
$\{E_i^{(m)}, F_i^{(m)}, K_{\mu} \, | 
\, i \in \tI (\rm{resp. } I), m \in \bbN, \mu \in Q\}$. 
We set 
\begin{eqnarray*}
 \left[
\begin{array}{r}
K_{\al_i} ; r\\
m  \, \, \,
\end{array}
\right]
:=\prod_{s=1}^m \frac{K_{\al_i}q^{r-s+1}-K_{\al_i}^{-1}q^{-r+s-1}}{q^s-q^{-s}},
\end{eqnarray*}
for $m \in \bbN, r \in \bbZ, i \in I$. 
It is known that 
$ \left[
\begin{array}{r}
K_{\al_i} ; r\\
m \, \,
\end{array}
\right] 
\in U_{\AA}^{\rm{res}}$
(see \cite{CP94b} \S 9.3A). 
By (\ref{pro RA}) and (\ref{pro RA-}), 
we have $T_i^{\pm 1}(\tU_{\AA}^{\rm{res}}) \subset\tU_{\AA}^{\rm{res}}$ 
(resp. $T_i^{\pm 1}(U_{\AA}^{\rm{res}}) \subset U_{\AA}^{\rm{res}}$)
for any $i \in \tI$ (resp. $I$). 
Hence we obtain 
$E_{\b}, F_{\b}, X_{i,r}^{\pm}, H_{i,s} \in \tU_{\AA}^{\rm{res}}$ 
for $\b \in \tD_{+}(I)$, $i \in I$, $r,s \in \bbZ(s \neq 0)$ 
by (\ref{def ARV}) and (\ref{rem DRV}). 
Similarly, we obtain $\bar{E}_{\ga}, 
\bar{F}_{\ga} \in U_{\AA}^{\rm{res}}$ 
for any $\ga \in \D_{+}$ by (\ref{def root vectors}). 
We define 
\begin{eqnarray*}
 \tU_{\e}^{\rm{res}}:=\tU_{\AA}^{\rm{res}} \otimes_{\AA} \bbC_{\e} 
\q  (\rm{resp. }U_{\e}^{\rm{res}}:=U_{\AA}^{\rm{res}} 
\otimes_{\AA} \bbC_{\e}). 
\end{eqnarray*}
Then we call $\tU_{\e}^{\rm{res}}$ (resp. $U_{\e}^{\rm{res}}$) 
``quantum loop algebra (resp. quantum algebra) of restricted type 
(or Lusztig type)'' (see (\cite{L89}, \cite{CP97})). 
We denote $E_{\b} \otimes 1$ by $e_{\b}$ for $\b \in \tD_{+}(I)$. 
Similarly, we set 
$f_{\b}:=F_{\b} \otimes 1$, $k_{\mu}:=K_{\mu} \otimes 1$, 
$x_{i,r}:=X_{i,r} \otimes 1$, $\cdots \cdots$. 
Moreover, we set 
\begin{eqnarray*}
&& \dot{e}_{\b}:=e_{\b}, 
\q \dot{f}_{\b}:=\Om (\dot{e}_{\b}), 
\q \dot{e}_{(i,m\d)}:=(\frac{1}{[m]}E_{(i,m\d)}) \otimes 1, 
\q \dot{f}_{(i,m\d)}:=\Om (\dot{e}_{(i,m\d)}),\\
&& \dot{h}_{i,s}:=(\frac{1}{[s]}H_{i,s}) \otimes 1, 
\q  \left[
\begin{array}{r}
k_{\al_i} ; r\\
m \, \,
\end{array}
\right] 
:= \left[
\begin{array}{r}
K_{\al_i} ; r\\
m \, \,
\end{array}
\right] \otimes 1, 
\end{eqnarray*} 
for $\b \in \tD_{+}^{\rm{re}}$, $i \in I$, $m \in \bbN$, 
$r,s \in \bbZ (s \neq 0)$.
Let $(\tU_{\e}^{\rm{res}})^{\pm}$ 
(resp. $(\tU_{\e}^{\rm{res}})^{0}$) 
be the $\bbC$-subalgebra of $\tU_{\e}^{\rm{res}}$ 
generated by 
$\{(x_{i,r}^{\pm})^{(m)} \, | \, i \in I, r \in \bbZ, m \in \bbN\}$ 
(resp. 
$\{k_{\mu}, 
 \left[
\begin{array}{r}
k_{\al_i} ; r\\
m \, \,
\end{array}
\right], 
\dot{h}_{i,s} \, 
|\, \mu \in Q, i \in I, r \in \bbZ, m \in \bbN, s \in \bbZ^{\times}\}$) 
and $(U_{\e}^{\rm{res}})^{+}$ 
(resp. $(U_{\e}^{\rm{res}})^{-}, (U_{\e}^{\rm{res}})^{0}$) 
be the $\bbC$-subalgebra of $U_{\e}^{\rm{res}}$ 
generated by $\{e_i^{(m)} \, | \, i \in I, m \in \bbN\}$ 
(resp. $\{f_i^{(m)} \, | \, i \in I, m \in \bbN\}$, 
$\{k_{\mu}, 
 \left[
\begin{array}{r}
k_{\al_i} ; r\\
m \, \,
\end{array}
\right] \, 
|\, \mu \in Q, i \in I, r \in \bbZ, m \in \bbN\}$). 
We obtain that $\left[
\begin{array}{r}
k_{\al_i} ; r\\
m \, \,
\end{array}
\right]$ 
is generated by $\{k_{\al_i}, 
\left[
\begin{array}{r}
k_{\al_i} ; 0\\
l \, \, \,
\end{array}
\right]\}$ for $i \in I, r \in \bbZ, m \in \bbN$.
It is known that 
$\Ures$ has the triangular decomposition, that is, 
the multiplication map defines an isomorphism of 
$\bbC$-vector spaces: 
\begin{eqnarray}
 (\Ures)^{-} \otimes (\Ures)^{0} \otimes (\Ures)^{+} \wt{\arr} \Ures 
\q (u^{-} \otimes u^0 \otimes u^{+} \mapsto u^{-}u^0u^{+}).
\label{triangular}
\end{eqnarray} 
Moreover, by \cite{CP97} Proposition 6.1, we have 
\begin{eqnarray}
 (\tUres)^{-} \otimes (\tUres)^{0} \otimes (\tUres)^{+}
\wt{\arr} \tUres 
\q (\wt{u}^{-} \otimes \wt{u}^0 \otimes \wt{u}^{+} 
\mapsto \wt{u}^{-}\wt{u}^0\wt{u}^{+}).
\label{affine triangular}
\end{eqnarray} 
\q We set 
\begin{eqnarray*}
&&(\tB_{\e}^{\rm{res}})^{+}
:=\{\prod_{\b \in \tD_{+}(I)}^{<} \dot{e}_{\b}^{(c(\b))}
 \, | \, c \in \bbZ^{\tD_{+}(I)}_{+}\}, 
\q (\tB_{\e}^{\rm{res}})^{-}:=\Om ((\tB_{\e}^{\rm{res}})^{+}), \\
&& (\tB_{\e}^{\rm{res}})^{0}
:=\{\prod_{i \in I}k_{\al_i}^{\d_i} 
 \left[
\begin{array}{r}
k_{\al_i} ; 0\\
m_i
\end{array}
\right]
\, | \, m_i \in \bbN, \d_i=0 \rm{ or } 1\}, 
\q \tB_{\e}^{\rm{res}} 
:=(\tB_{\e}^{\rm{res}} )^{-}(\tB_{\e}^{\rm{res}} )^{0}
(\tB_{\e}^{\rm{res}} )^{+},
\end{eqnarray*}
(see (\ref{sig APBW})).
By \cite{G}, we obtain the following theorem. 
\newtheorem{thm RSAPBW}[sec 3]{Theorem}
\begin{thm RSAPBW}
\label{thm RSAPBW}
$\tB_{\e}^{\rm{res}}$ is a $\bbC$-basis of $\tU_{\e}^{\rm{res}}$.
\end{thm RSAPBW}
Proof. 
By \cite{G}, we obtain a PBW basis of $\tU_{\AA}^{\rm{res}}$.
Since any $\AA$-basis of $\tU_{\AA}^{\rm{res}}$ become 
$\bbC$-basis of $\tU_{\AA}^{\rm{res}} \otimes _{\AA} \bbC_{\e}$ 
canonically (see Lemma \ref{lem basis}), 
we obtain this theorem. 
\qed
\subsection{Small quantum algebras}
\q Let $\tUfin$ (resp. $\Ufin$) 
be the $\bbC$-subalgebra of $\tUres$ (resp. $\Ures$) 
generated by 
$\{e_i, f_i, k_{\mu} \, | \, i \in \tI (\rm{ resp. } I), \mu \in Q\}$). 
Then we call $\tUfin$ (resp. $\Ufin$) 
``small quantum loop algebra (resp. small quantum algebra)''. 
Let $\DD_{\e}^{\rm{fin}}$ be the $\bbC$-subalgebra of $\tUres$ 
generated by 
$\{x_{i,r}^{\pm}, h_{i,s}, k_{\mu} \, | \, 
i \in I, r,s \in \bbZ(s \neq 0), \mu \in Q\}$. 
By (\ref{pro DR}), we obtain $e_0, f_0 \in \DD_{\e}^{\rm{fin}}$. 
So we have $\DD_{\e}^{\rm{fin}}=\tUfin$. 
Let $(\tUfin)^{\pm}$ (resp. $(\tUfin)^{0}$) 
be the $\bbC$-subalgebra of $\tUfin$ 
generated by $\{x_{i,r}^{\pm} \, | \, i \in I, r \in \bbZ\}$ 
(resp. 
$\{h_{i,s}, k_{\mu} \, | \, i \in I, s \in \bbZ^{\times}, \mu \in Q\}$) 
and $(\Ufin)^{+}$ (resp. $(\Ufin)^{-}, (\Ufin)^{0}$) 
be the $\bbC$-subalgebra of $\Ufin$ 
generated by 
$\{e_i\}_{i \in I}$ (resp. $\{f_i\}_{i \in I}, 
\{k_{\mu}\}_{\mu \in Q}$). 
Let $\bbZ_l^{\D_{+}}$ be as in (\ref{def blpbw}).
We set 
\begin{eqnarray}
&& (\Bfin)^{+}:=\{\prod_{\ga \in \D_{+}}^{<} \bar{e}_{\ga}^{c(\ga)} \,
| \, c \in \bbZ_l^{\D_{+}}\}, 
\q (\Bfin)^{-}:=\Om((\Bfin)^{+}), \no \\ 
&& (\Bfin)^{0}:=\{k_{\mu} \, | \, \mu \in Q_{2l}\}, 
\q \Bfin:=(\Bfin)^{-}(\Bfin)^{0}(\Bfin)^{+}. 
\label{def bfin}
\end{eqnarray}  
\newtheorem{thm SPBW}[sec 3]{Theorem}
\begin{thm SPBW}[\cite{L90a} \S5, \cite{L90b} \S8]
\label{thm SPBW}
$(\Bfin)^{\star}$ (resp. $\Bfin$) 
is a $\bbC$-basis of $(\Ufin)^{\star}$ (resp. $\Ufin$) for $\s=$.
\end{thm SPBW}
Since $e_{\b}, f_{\b}$ are generated by 
$\{e_i\}_{i \in \tI}, \{f_i\}_{i \in \tI}$ respectively, 
we get $e_{\b}, f_{\b} \in \tUfin$ for any $\b \in \tD_{+}(I)$. 
For $\b \in \tD_{+}^{\rm{re}}$, $i \in I$, $m \in \bbN$, 
we have 
\begin{eqnarray*}
&& e_{\b}^l=[l]_{\e}! \dot{e}_{\b}^{(l)}, 
\q f_{\b}^l=[l]_{\e}! \dot{f}_{\b}^{(l)}, 
\q e_{(i,ml\d)}=[l]_{\e}\dot{e}_{(i,ml\d)}, 
\q f_{(i,ml\d)}=[l]_{\e}\dot{f}_{(i,ml\d)}, \\ 
&& \prod_{r=1}^l(k_{\al_i} \e^{1-r}-k_{\al_i} \e^{r-1})
=\prod_{r=1}^l (\e^r-\e^{-r})
 \left[
\begin{array}{r}
k_{\al_i} ; 0\\
l \, \,
\end{array}
\right],
\end{eqnarray*} 
(see \cite{L89} Lemma 4.4). 
So we obtain 
\begin{eqnarray}
 e_{\b}^l=f_{\b}^l=e_{(i,sl\d)}=f_{(i,sl\d)}=k_{\al_i}^{2l}-1=0, 
\label{pro NE}
\end{eqnarray}
for any $\b \in \tD_{+}^{\rm{re}}$, $s \in \bbN$, $i \in I$. 
Let $\bbZ_l^{\tD_{+}(I)}$ be as in (\ref{def ARSl-basis}).
Define
\begin{eqnarray*}
&& (\tBfin)^{+}:=\{\prod_{\b \in \tD_{+}(I)}^{<} e_{\b}^{c(\b)} \, 
| \, c \in \bbZ_l^{\tD_{+}(I)}\}, 
\q (\tBfin)^{-}:=\Om((\tBfin)^{+}), \\
&& (\tBfin)^{0}:=\{k_{\mu} \, | \, \mu \in Q_{2l}\}, 
\q \tBfin:=(\tBfin)^{-}(\tBfin)^{0}(\tBfin)^{+}.
\end{eqnarray*} 
Since $e_{\b}, f _{\b} \in \tUfin$ for any $\b \in \tD_{+}(I)$, 
we have $\tBfin \subset \tUfin$. 
Therefore, by Theorem \ref{thm RSAPBW}, 
we obtain the following lemma.
\newtheorem{pro ASPBW}[sec 3]{Lemma}
\begin{pro ASPBW}
\label{pro ASPBW}
$\tBfin$ is linearly independent in $\tUfin$.
\end{pro ASPBW}
\subsection{Isomorphism theorem}
\newtheorem{pro ISO}[sec 3]{Proposition}
\begin{pro ISO}[\cite{AN} Lemma 4.8]
\label{pro ISO}
There exists a $\bbC$-algebra isomorphism 
$\bar{\phi}: U_{\e}/I_{\e} \arr \Ufin$ 
such that $\bar{\phi}(E_i+I_{\e})=e_i$, $\bar{\phi}(F_i+I_{\e})=f_i$ 
and $\bar{\phi}(K_{\mu}+I_{\e})=k_{\mu}$ 
for $i \in I, \mu \in Q$. 
\end{pro ISO}
We obtain an affine version of the above result:
\newtheorem{thm ISO}[sec 3]{Theorem}
\begin{thm ISO}
\label{thm ISO}
There exists a $\bbC$-algebra isomorphism 
$\widetilde{\phi}: \tU_{\e}/\tI_{\e} \arr \tUfin$ 
such that $\widetilde{\phi}(E_i+\tI_{\e})=e_i$, 
$\widetilde{\phi}(F_i+\tI_{\e})=f_i$ 
and $\widetilde{\phi}(K_{\mu}+\tI_{\e})=k_{\mu}$ 
for $i \in \tI, \mu \in Q$. 
In particular, $\tBfin$ is a $\bbC$-basis of $\tUfin$.
\end{thm ISO}
Proof. Elements in $\{e_i, f_i, k_{\mu} \, | \, 
i \in \tI, \mu \in Q\} (\subset \tUfin)$ 
satisfy the relations of Definition \ref{def QA}. 
So, by the universality of $\tU_{\e}$ 
(see Remark \ref{rem NS} (b)), 
there exists a surjective $\bbC$-algebra homomorphism 
$\phi: \tU_{\e} \arr \tUfin$ 
such that $E_i \mapsto e_i, F_i \mapsto f_i, K_{\mu} \mapsto k_{\mu}$ 
for $i \in \tI, \mu \in Q$. 
Since $e_{\b}=E_{\b} \otimes 1$ and $f_{\b}=F_{\b} \otimes 1$ 
($\b \in \tD_{+}(I)$), 
we obtain $\phi(E_{\b})=e_{\b}$ and $\phi(F_{\b})=f_{\b}$. 
Then, by (\ref{pro NE}), we have $\phi(\tI_{\e})=0$. 
Hence there exists a surjective $\bbC$-algebra homomorphism 
$\widetilde{\phi}: \tue / \tI_{\e} \arr \tUfin$ 
such that $\widetilde{\phi}(E_{\b}+\tI_{\e})=e_{\b}$, 
$\widetilde{\phi}(F_{\b}+\tI_{\e})=f_{\b}$ and 
$\widetilde{\phi}(K_{\mu}+\tI_{\e})=k_{\mu}$ 
for $\b \in \tD_{+}(I)$, $\mu \in Q$. \\
\q Let $u \in \rm{Ker}(\wt{\phi})$. 
By Proposition \ref{pro lAPBW}, 
we have 
\begin{eqnarray*}
 u=\sum_{\mu \in Q_{2l}}\sum_{c, c^{'} \in \bbZ_l^{\tD_{+}(I)}}
a(c, \mu,c^{'})(\prod_{\b \in \tD_{+}(I)}^{<}E_{\b}^{c(\b)}+\tI_{\e}) 
(K_{\mu}+\tI_{\e})
(\prod_{\b \in \tD_{+}(I)}^{<}F_{\b}^{c^{'}(\b)}+\tI_{\e}),
\end{eqnarray*}
where $a(c, \mu, c^{'}) \in \bbC$. Then we get
\begin{eqnarray*}
0= \widetilde{\phi}(u)=
\sum_{\mu \in Q_{2l}}\sum_{c, c^{'} \in \bbZ_l^{\tD_{+}(I)}}
a(c, \mu,c^{'})(\prod_{\b \in \tD_{+}(I)}^{<}e_{\b}^{c(\b)}) 
k_{\mu}(\prod_{\b \in \tD_{+}(I)}^{<}f_{\b}^{c^{'}(\b)}).  
\end{eqnarray*}
Hence, by Lemma \ref{pro ASPBW}, we obtain
$a(c,\mu, c^{'})=0$ for any 
$c,c^{'} \in \bbZ_l^{\tD_{+}(I)}, \mu \in
Q_{2l}$. 
Thus $\widetilde{\phi}$ is injective. 
Therefore $\widetilde{\phi}$ is an isomorphism 
and $\tBfin$ is a $\bbC$-basis of $\tUfin$. 
\qed 
\section{Evaluation representations of restricted type}
\setcounter{equation}{0}
\renewcommand{\theequation}{\thesection.\arabic{equation}}
\subsection{Representation theory of restricted type} 
\q We call a $\tUres$-module $\tV$ (resp. $\Ures$-module $V$) 
``type \bf{1}'' if $k_{\mu}^l=1$ on $\tV$ (resp. $V$) 
for any $\mu \in Q$. 
In general, finite dimensional irreducible $\tUres$-modules  
(resp. $U_{\e}^{\rm{res}}$-modules) are classified
into $2^n$ types 
according to $\{\sigma : Q \arr \{\pm 1\}; \rm{ homomorphism of group
}\}$. 
It is known that for any $\sigma : Q \arr \{\pm 1\}$, 
the category of finite dimensional $\tUres$-module 
(resp. $\Ures$-module) of type $\sigma$ is essentially equivalent to 
the category of the finite dimensional $\tUres$-module (resp. $\Ures$-module) 
of type $\bf{1}$. 

Let $U=\Ures, \Ufin, \tUres$ or $\tUfin$. 
\newtheorem{def sec4}{Definition}[section]
\begin{def sec4} 
\label{def sec4}
Let $V$ be a $U$-module and $v$ be a nonzero vector in $V$. 
Suppose that $v$ satisfies 
\begin{eqnarray*}
&&\hspace{-30pt}e_i^{(m)}v=0  \text{ for any }i\in I,\,\,m\in\bbN
\,\, \rm{if }  U=\Ures, 
\,\, e_iv=0\text{ for any }i\in I\text{ if }  U=\Ufin, \\
&&\hspace{-30pt}
(x_{i,r}^{+})^{(m)}v=0
\text{ for any }i\in I,\,\,m\in\bbN,\,\,r\in\bbZ \,\,  \rm{if } U=\tUres, 
\,\, x_{i,r}^{+}v=0
\text{ for any }i\in I,\,\,r\in\bbZ
\,\, \rm{if } U=\tUfin, 
\end{eqnarray*}
We call $v$ a ``primitive vector'' in $V$.
\end{def sec4}
\newtheorem{def RSHW}[def sec4]{Definition}
\begin{def RSHW} 
\label{def RSHW}
Let $V$ be a $U$-module  
and $\L: U^0 \arr \bbC$ be a $\bbC$-algebra homomorphism.
We assume that $V$ is generated as a $U$-module by
a primitive vector $v_{\L} \in V$ such that 
\begin{eqnarray*}
u_0v_{\L}=\L(u_0)v_{\L},
\end{eqnarray*}
for any $u_0 \in U^0$. 
Then we call $V$ a ``highest weight $U$-module''  
generated by a ``highest weight vector'' $v_{\L}$ 
with ``highest weight $\L$''. 
\end{def RSHW}
\newtheorem{pro HWM}[def sec4]{Proposition}
\begin{pro HWM} 
\label{pro HWM}
For any $\bbC$-algebra homomorphism $\L: U^0 \arr \bbC$, 
there exists a unique (up to isomorphism) 
irreducible highest weight $U$-module $V$ with highest
 weight $\L$. 
\end{pro HWM}
Proof. 
For any $\Ures$ (resp. $\tUres$-module) $V$, 
we can define the  weight spaces on $V$ by 
\begin{eqnarray}
 V_{\mu}:=\{v \in V \, | \, k_{\al_i}v=\e^{ \langle\mu, 
\al_i^{\vee} \rangle}v, 
\q \left[
\begin{array}{r}
k_{\al_i} ; 0\\
l \, \, \,
\end{array}
\right]v
= \left[
\begin{array}{r}
\langle \mu, \al_i^{\vee} \rangle\\
l \, \, \q
\end{array}
\right]_{\e}v 
\q \rm{for any $i \in I$}\}, 
\label{def weight space2}
\end{eqnarray}
where $\mu \in P$ (see \cite{L89}, \cite{CP97}). 
Then, by the theory of highest weight modules, 
we obtain this proposition in the case of $U=\Ures$ or $\tUres$ 
(see \cite{CP97} Proposition 7.3). 
So we shall prove the case of $\tUfin$. 
We can prove the $\Ufin$ case similarly. 

Let $\hUfin$ be the $\bbC$-subalgebra of $\tUres$ 
generated by $\tUfin \cup 
\{\left[
\begin{array}{r}
k_{\al_i} ; 0\\
l \, \, \,
\end{array}
\right] \, | \, i \in I\}$. 
For any $\bbC$-algebra homomorphism $\L: (\tUfin)^0 \arr \bbC$, 
let  $\wh{I}_{\e}^{\rm{fin}}(\L)$ (resp. $\tI_{\e}^{\rm{fin}}(\L)$) 
be the left ideal of $\hUfin$ (resp. $\tUfin$) generated by 
$\{x_{i,r}^{+}, u_0-\L(u_0),
 \left[
\begin{array}{r}
k_{\al_i} ; 0\\
l \, \, \,
\end{array}
\right]
\, | \, i \in I, r \in \bbZ, u_0 \in (\tUfin)^{0} \}$ 
(resp. $\{x_{i,r}^{+}, u_0-\L(u_0)  
\, | \, i \in I, r \in \bbZ, u_0 \in (\tUfin)^{0} \}$).
We define a $\hUfin$-module $\hMfin(\L)$ 
and a $\tUfin$-module $\tMfin(\L)$ respectively by 
\begin{eqnarray*}
 \hMfin(\L):=\hUfin /  \hIfin(\L),
\q  \tMfin(\L):=\tUfin / \tIfin(\L).
\end{eqnarray*} 
We set $\wh{v}_{\L}:=1+\hIfin(\L) \in \hMfin(\L)$ 
and $\wt{v}_{\L}:=1+\tIfin(\L) \in \tMfin(\L)$. 
Let $\hNfin(\L)$ be the $\tUfin$-submodule of $\hMfin(\L)$ 
generated by $\wh{v}_{\L}$.  
Then, by the universality of the $\tUfin$-module $\tMfin(\L)$, 
there exists a surjective $\tUfin$-module homomorphism 
$\phi: \tMfin(\L) \arr \hNfin(\L)$ 
such that $\phi(\wt{v}_{\L})=\wh{v}_{\L}$.
Let $B$ be a $\bbC$-basis of $(\tUfin)^{-}$. 
Then, by (\ref{affine triangular}), 
we obtain that $\{u\wh{v}_{\L} \, | \, u \in B\}$ 
(resp. $\{u\wt{v}_{\L} \, | \, u \in B\}$) 
is a $\bbC$-basis of $\hMfin(\L)$ (resp. $\tMfin(\L)$). 
Hence $\phi$ is an isomorphism of $\tUfin$-module. 
So we can regard $\tMfin(\L)$ as $\hUfin$-module. 
By the similar way to (\ref{def weight space2}), 
we can define the  weight spaces on this module. 
Then, by the theory of the  highest weight module, 
$\tMfin(\L)$ has a unique simple quotient of $\tUfin$-module 
and it is the unique irreducible highest weight $\tUfin$-module 
with highest weight $\L$. 
\qed 

For any $\bbC$-algebra homomorphism $\L : U^0 \arr \bbC$, 
we denote the unique irreducible highest weight $U$-module 
with highest weight $\L$ by 
$V_{\e}^{\rm{res}}(\L)$ if $U=\Ures$, 
$V_{\e}^{\rm{fin}}(\L)$ if $U=\Ufin$, 
$\wt{V}_{\e}^{\rm{res}}(\L)$ if $U=\tUres$, 
and $\wt{V}_{\e}^{\rm{fin}}(\L)$ if $U=\tUfin$. 
Then, by Proposition \ref{pro HWM} 
and the uniqueness of the primitive vectors,
we obtain the following proposition. 
\newtheorem{pro classification of HWM}[def sec4]{Proposition}
\begin{pro classification of HWM}
\label{pro classification of HWM}  
Let $\L$ and $\L^{'}: U^0 \arr \bbC$ 
be $\bbC$-algebra homomorphisms. 
Then $V_{\e}^{\rm{res}}(\L)$ 
(resp. $V_{\e}^{\rm{fin}}(\L)$, $\wt{V}_{\e}^{\rm{res}}(\L)$, 
$\wt{V}_{\e}^{\rm{fin}}(\L)$) 
is isomorphic to $V_{\e}^{\rm{res}}(\L^{'})$ 
(resp. $V_{\e}^{\rm{fin}}(\L^{'})$, $\wt{V}_{\e}^{\rm{res}}(\L^{'})$, 
$\wt{V}_{\e}^{\rm{fin}}(\L^{'})$) 
if and only if $\L=\L^{'}$.   
\end{pro classification of HWM}
Now, we define $\PP_{i,m} \in \tU_q$  
inductively by 
\begin{eqnarray*}
\PP_{i,0}:=1, 
\q \PP_{i,m}:=-\frac{K_{\al_i}^{-1}}{1-q^{-2m}}
\sum_{s=1}^m \Psi_{i,s}^{+}\PP_{i,m-s}, 
\q \PP_{i,-m}:=\Omega(\PP_{i,m}),
\end{eqnarray*}
where $\Omega$ as in (\ref{pro omega}).
We have $\Om(\Psi_{i,s}^{+})=\Psi_{i,-s}^{-}$. 
\newtheorem{pro PP}[def sec4]{Proposition}
\begin{pro PP}[\cite{CP97} \S 3]
\label{pro PP}
For any $i \in I$, $r \in \bbZ$, 
we have $\PP_{i,r} \in \tU_{\AA}^{\rm{res}}$. 
Moreover, $(\tUres)^0$ is generated by 
$\{k_{\al_i}, 
\left[
\begin{array}{r}
k_{\al_i} ; 0\\
l \, \, \,
\end{array}
\right], 
\PP_{i,r}\otimes 1 \, | \, i \in I, r \in \bbZ \}$ as $\bbC$-algebra. 
\end{pro PP}
We simply denote $\PP_{i,r} \otimes 1 \in \tUres$ by $\PP_{i,r}$.
We set 
\begin{eqnarray*}
 \bbC_0[t]:=\{P \in \bbC[t] \, | \, 
P \rm{ is monic }, P(0) \neq 0\}.
\end{eqnarray*}
We call a polynomial $P \in \bbC[t]$ ``$l$-acyclic'' 
if it is not divisible by $(1-ct^l)$ for any $c \in \bbC^{\times}$ 
(see \cite{F}) and set 
\begin{eqnarray*}
 \bbC_l[t]:=\{P \in \bbC_0[t] \, | \, \rm{ $P$ is $l$-acyclic }\}.
\end{eqnarray*}
\newtheorem{def fin dim HW}[def sec4]{Definition}
\begin{def fin dim HW}
\label{def fin dim HW}
(a) For $\l=(\l_i)_{i \in I} \in \bbZ_{+}^n$, 
let $\l_i^{(0)} \in \bbZ_l$ and $\l_i^{(1)} \in \bbZ_{+}$ 
such that $\l_i=\l_i^{(0)}+l\l_i^{(1)}$ ($i \in I$). 
We define a $\bbC$-algebra homomorphism 
$\L_{\l}^{\rm{res}}: (\Ures)^0 \arr \bbC$ by 
\begin{eqnarray*}
&& \L_{\l}^{\rm{res}}(k_{\al_i}):=\e^{\l_i^{(0)}}, 
\q \L_{\l}^{\rm{res}}(
\left[
\begin{array}{r}
k_{\al_i} ; 0\\
l \q
\end{array}
\right])
:=\l_i^{(1)} 
\q (i \in I).
\end{eqnarray*}
\q (b) For $\l=(\l_i)_{i \in I} \in \bbZ_{l}^n$, 
we define a $\bbC$-algebra homomorphism 
$\L_{\l}^{\rm{fin}}: (\Ufin)^0 \arr \bbC$ by 
\begin{eqnarray*}
&& \L_{\l}^{\rm{fin}}(k_{\al_i}):=\e^{\l_i} 
\q (i \in I). 
\end{eqnarray*}
\q (c) For $\bf{P}=(P_i)_{i \in I} \in \bbC_0[t]^n$, 
let $p_i^{(0)} \in \bbZ_l$ and $p_i^{(1)} \in \bbZ_{+}$ 
such that $\rm{deg}(P_i)=p_i^{(0)}+lp_i^{(1)}$ ($i \in I$). 
We define a $\bbC$-algebra homomorphism 
$\wt{\L}_{\bf{P}}^{\rm{res}}: (\tUres)^0 \arr \bbC$ by 
\begin{eqnarray*}
&& \wt{\L}_{\bf{P}}^{\rm{res}}(k_{\al_i}):=\e^{p^{(0)}_i}, 
\q \wt{\L}_{\bf{P}}^{\rm{res}}(
\left[
\begin{array}{r}
k_{\al_i} ; 0\\
l \q
\end{array}
\right])
:=p_i^{(1)} 
\q (i \in I), \\
&&\sum_{m=0}^{\infty}\wt{\L}_{\bf{P}}^{\rm{res}}
(\PP_{i,m})t^m:=\frac{P_i(t)}{P_i(0)}, 
\q \sum_{m=0}^{\infty}\wt{\L}_{\bf{P}}^{\rm{res}}
(\PP_{i,-m})t^m:=\frac{Q_i(t)}{Q_i(0)}, 
 \end{eqnarray*}
where $Q_i(t):=t^{\rm{deg}(P_i)}P_i(t^{-1})$ 
(see \cite{CP97} \S8).  \\
\q (d) For $\bf{P}=(P_i)_{i \in I} \in \bbC_l[t]^n$, 
we define a $\bbC$-algebra homomorphism 
$\wt{\L}_{\bf{P}}^{\rm{fin}}: (\tUfin)^0 \arr \bbC$ by 
\begin{eqnarray*}
\sum_{m=0}^{\infty}\wt{\L}_{\bf{P}}^{\rm{fin}}
(\psi_{i,m}^{+})t^m 
:=\e^{\rm{deg}(P_i)}\frac{P_i(\e^{-2}t)}{P_i(t)}
:= \sum_{m=0}^{\infty}\wt{\L}_{\bf{P}}^{\rm{fin}}
(\psi_{i,-m}^{-})t^{-m}, 
\label{fac DP} 
\end{eqnarray*}
in the sense that the left and right hand sides are 
the Laurent expansions
of the middle term about $t=0$ and $t=\infty$, respectively 
(see \cite{CP97} \S8). 
\end{def fin dim HW}
By \cite{CP97} \S8, 
we obtain $\wt{\L}_{\bf{P}}^{\rm{fin}}
=\wt{\L}_{\bf{P}}^{\rm{res}}|_{(\tUfin)^0}$ ($\bf{P} \in \bbC_l[t]^n$).
For $\l \in \bbZ_{+}^{n}$ (resp. $\bbZ_l^n$), 
we set $\Vres:=V_{\e}^{\rm{res}}(\L^{\rm{res}}_{\l})$ 
(resp. $\Vfin:=V_{\e}^{\rm{fin}}(\L^{\rm{fin}}_{\l})$). 
Similarly, for $\bf{P} \in \bbC_0[t]^n$ (resp. $\bbC_l[t]^n$), 
we set 
$\tVres:=\wt{V}_{\e}^{\res}(\wt{\L}^{\rm{res}}_{\bf{P}})$ 
(resp. 
$\tVfin:=\wt{V}_{\e}^{\fin}(\wt{\L}^{\rm{fin}}_{\bf{P}})$) and
call $\bf{P}$ ``Drinfel'd polynomial'' of $\tVres$ (resp. $\tVfin$).
\newtheorem{thm RSCT}[def sec4]{Theorem}
\begin{thm RSCT}[\cite{L89}, \cite{CP94b}Proposition 11.2.10]
\label{thm RSCT}
For any $\l \in \bbZ_{+}^n$ (resp. $\bbZ_l^{+}$), 
$\Vres$ (resp. $\Vfin$) is a finite dimensional irreducible 
$\Ures$-module (resp. $\Ufin$-module) of type \bf{1}. 
Conversely, for any finite dimensional irreducible 
$\Ures$-module (resp. $\Ufin$-module) $V$ of type \bf{1}, 
there exists a unique $\l \in \bbZ_{+}^n$ (resp. $\bbZ_l^n$)
such that $V$ is isomorphic to $\Vres$ (resp. $\Vfin$) 
as a $\Ures$-module (resp. $\Ufin$-module). 
In particular, for any $\l \in \bbZ_l^n$, 
$\Vfin$ is isomorphic to $\Vres$ as a $\Ufin$-module.
\end{thm RSCT}
\newtheorem{thm ARSCT}[def sec4]{Theorem}
\begin{thm ARSCT}[\cite{CP97} Theorem 8.2, 9.2, \cite{F} Theorem 2.6]
\label{thm ARSCT}
For any $\bf{P} \in \bbC_{0}[t]^n$ (resp. $\bbC_l[t]^n$), 
$\tVres$ (resp. $\tVfin$) is a finite dimensional irreducible 
$\tUres$-module (resp. $\tUfin$-module) of type \bf{1}. 
Conversely, for any finite dimensional irreducible 
$\tUres$-module (resp. $\tUfin$-module) $V$ of type \bf{1}, 
there exists a unique $\bf{P} \in \bbC_{0}[t]^n$ (resp. $\bbC_l[t]^n$)
such that $V$ is isomorphic to $\tVres$ (resp. $\tVfin$) 
as $\tUres$-module (resp. $\tUfin$-module). 
In particular, for any $\bf{P} \in \bbC_l[t]^n$, 
$\tVfin$ is isomorphic to $\tVres$ as $\tUfin$-module.
\end{thm ARSCT}
By the tensor product theorem, in order to understand 
the representation theory of $\tUres$ (resp. $\Ures$),
it is sufficient to consider
 $\tUfin$ (resp. $\Ufin$) 
(see \cite{CP97} and \cite{L89}). 
\subsection{Drinfel'd polynomials of evaluation representations}
\q For $m \in \bbZ_{+}$, 
let $V_{\e}^{\rm{fin}}(m)$ be the $(m+1)$-dimensional irreducible 
$\Ufin(\sl_2)$-module. 
By (\ref{pro embed}), we can regard $\tVfin$ as $\Ufin(\sl_{n+1})$-module 
($\bf{P} \in \bbC_l[t]^n$). 
Then, by \cite{CP97} \S 7-9 (in particular, p.321), 
we obtain the following theorem. 
\newtheorem{thm sl2CT}[def sec4]{Theorem}
\begin{thm sl2CT}[\cite{CP97}]
\label{thm sl2CT}
For any $P \in \bbC_l[t] (P \neq 1)$, 
there exist $r, m_s \in \bbN$ and
 $c_s \in \bbC^{\times} (1 \leq s \leq r)$ 
such that 
\begin{eqnarray*}
 \tV^{\rm{fin}}_{\e}(P) \cong \tV^{\rm{fin}}_{\e}(P_1)  \otimes \cdots 
\otimes \tV^{\rm{fin}}_{\e}(P_r)  
\qq \rm{(as a $\Ufin(\wt{\sl}_2)$-module) },
\end{eqnarray*}
where $P=\prod_{s=1}^rP_s$,
$P_s(t)=\prod_{p=1}^{m_s}(t-c_s\e^{m_s+1-2p}) \in \bbC_l[t]$ 
($1 \leq s \leq r$). 
In particular, $\tV_{\e}^{\rm{fin}}(P_s)$ is isomorphic to 
$V_{\e}^{\rm{fin}}(m_s)$ as a $\Ufin(\sl_2)$-module. 
\end{thm sl2CT}
By this theorem, we obtain the following lemma.
\newtheorem{lem slSR}[def sec4]{Lemma}
\begin{lem slSR} 
\label{lem slSR}
Let $P \in \bbC_l[t]$, $m \in \bbN$. 
If $\tV^{\rm{fin}}_{\e}(P)$ is isomorphic to $V_{\e}^{\rm{fin}}(m)$ 
as a $\Ufin(\sl_2)$-module, 
then there exists $c \in \bbC^{\times}$ 
such that $P(t)=\prod_{p=1}^{m}(t-c\e^{m+1-2p})$. 
\end{lem slSR}
By using this lemma, we can prove the following lemma.
\newtheorem{lem ER}[def sec4]{Lemma}
\begin{lem ER} 
\label{lem ER}
Let $\bf{P}=(P_i)_{i \in I} \in \bbC_l[t]^n$, 
$\l=(\l_i)_{i \in I} \in \bbZ_l^n$. 
We assume $\tVfin$ is isomorphic to $\Vfin$ as a $\Ufin$-module. 
Then, for any $i \in I (P_i \neq 1)$, 
there exists $c_i \in \bbC^{\times}$ such that 
$P_i(t)=\prod_{p=1}^{\l_i}(t-c_i\e^{\l_i+1-2p})$. 
\end{lem ER}
Proof. 
Let $v_{\bf{P}}$ be a highest weight vector in $\tVfin$. 
Then $v_{\bf{P}}$ is also a highest weight vector in $\Vfin$. 
For $i \in I$, 
let $(\tUfin)_i$ (resp. $(\Ufin)_i$) be the $\bbC$-subalgebra of 
$\tUfin$ (resp. $\Ufin$) generated by 
$\{x_{i,r}^{\pm}, h_{i,s}, k_{\al_i} \, | 
r,s \in \bbZ, s \neq 0\}$ (resp. $\{e_i, f_i, k_{\al_i}\}$) 
and $\widetilde{W}_i$ (resp. $W_i$) be the $(\tUfin)_i$-submodule 
(resp. $(\Ufin)_i$-submodule) of $\tVfin$ generated by $v_{\bf{P}}$. 
By (\ref{pro embed2}) (resp. (\ref{pro embed})), 
we can regard $\widetilde{W}_i$ (resp. $W_i$) 
as a $\Ufin(\widetilde{\sl}_2)$-module (resp. $\Ufin(\sl_2)$-module). 
Then, by Lemma 7.6 in \cite{CP97}, we obtain 
$\wt{W}_i \cong \tV_{\e}^{\rm{fin}}(P_i)$ as a $\Ufin(\wt{\sl}_2)$-module. 
Similarly (more easily), we obtain 
$W_i \cong V_{\e}^{\rm{fin}}(\l_i)$ as a $\Ufin(\sl_2)$-module. 
So, by Lemma \ref{lem slSR}, 
it is enough to prove $\wt{W}_i=W_i$ 
for $i \in I$ such that $P_i \neq 1$. 

By (\ref{affine triangular}),  
$\widetilde{W}_i$ is spanned by 
$\{x_{i,r_1}^{-} \cdots x_{i,r_m}^{-}v_{\bf{P}} \, | \, 
r_1, \cdots r_m \in \bbZ, m \in \bbZ_{+}\}$ as a $\bbC$-vector space. 
By Theorem \ref{thm ARSCT}, we can regard $\tVfin$ 
as a $\tUres$-module and define the  weight spaces on $\tVfin$ 
by the similar way to (\ref{def weight space2}). 
Then, by the relations of Drinfel'd realization, 
we have $x_{i,r_1}^{-} \cdots x_{i,r_m}^{-}v_{\bf{P}} 
\in (\tVfin)_{\l-m\al_i}$ for any $m \in \bbN, r_1, \cdots r_m \in
\bbZ$. 
So we obtain $\widetilde{W}_i \subset \bigoplus_{m \geq 0}
(\tVfin)_{\l-m\al_i}$. 
On the other hand, by the assumption of this lemma, 
$\tVfin$ is isomorphic to $\Vfin$ as a $\Ufin$-module. 
Hence, by (\ref{triangular}), 
$\tVfin$ is spanned by 
$\{\prod_{\ga \in \D_{+}}^{<}\bar{f}_{\ga}v_{\bf{P}} \, |
\, c \in \bbZ_l^{\D_{+}}\}$. 
Since $\prod_{\ga \in \D_{+}}^{<}\bar{f}_{\ga}^{c(\ga)}v_{\bf{P}} 
\in (\tVfin)_{\sum_{\ga \in \D_{+}}c(\ga)\ga}$, 
we have 
\begin{eqnarray*}
 \bigoplus_{m \geq 0}(\tVfin)_{\l-m\al_i} = 
\bigoplus_{m \in \bbZ_l} \bbC f_i^mv_{\bf{P}} \subset W_i. 
\end{eqnarray*}
Then we have $\widetilde{W}_i=W_i$. 
\qed   \\ \\
\q Let $(U_{\AA}^{\rm{res}})^{'}$ (resp. $(\Ufin)^{'}$) 
be the extended algebra of $U_{\AA}^{\rm{res}}$ (resp. $\Ufin$) 
defined by replacing $\{K_{\mu} \, | \, 
\mu \in Q\}$ with $\{K_{\mu} \, | \, \mu \in P\}$ 
(see Definition \ref{def QA}). 
By (\ref{pro EH2}), we obtain 
$\rm{ev}_{\bf{a}}^{\pm}(\tU_{\AA}^{\rm{res}}) 
\subset (U_{\AA}^{\rm{res}})^{'}$ ($\bf{a} \in \bbC^{\times}$).
Hence, by Proposition \ref{pro EH}, 
we obtain the evaluation $\tUfin$-homomorphisms 
$(\rm{ev}_{\bf{a}}^{\rm{fin}})^{\pm}: \tUfin \arr (\Ufin)^{'}$ 
defined by 
\begin{eqnarray}
&& (\rm{ev}_{\bf{a}}^{\rm{fin}})^{\pm}(e_i):=e_i, 
\q ( \rm{ev}_{\bf{a}}^{\rm{fin}})^{\pm}(f_i):=f_i 
\q (\rm{ev}_{\bf{a}}^{\rm{fin}})^{\pm}(k_{\mu}):=k_{\mu},  \no\\
&& (\rm{ev}_{\bf{a}}^{\rm{fin}})^{+}(e_0):=
\bf{a}\e^{-1}k_{\L_1}k_{\L_n}^{-1}[f_n,[f_{n-1}, \cdots, [f_2,f_1]_{\e^{-1}}
\cdots ]_{\e^{-1}}, \no \\
&& (\rm{ev}_{\bf{a}}^{\rm{fin}})^{-}(f_0):=
{\bf{a}}^{-1}(-1)^{n-1}\e^{n}k_{\L_1}^{-1}k_{\L_n}
[e_n,[e_{n-1}, \cdots, [e_2,e_1]_{\e^{-1}}
\cdots ]_{\e^{-1}},\no \\
&&(\rm{ev}_{\bf{a}}^{\rm{fin}})^{+}(e_0):=
\bf{a}\e^{-1}k_{\L_1}^{-1}k_{\L_n}[f_1,[f_{2}, \cdots, [f_{n-1},f_n]_{\e^{-1}}
\cdots ]_{\e^{-1}}, \no \\
&&(\rm{ev}_{\bf{a}}^{\rm{fin}})^{-}(f_0):=
{\bf{a}}^{-1}(-1)^{n-1}\e^{n}k_{\L_1}k_{\L_n}^{-1}
[e_1,[e_{2}, \cdots, [e_{n-1},e_n]_{\e^{-1}}
\cdots ]_{\e^{-1}},  
\label{sig SQEH}
\end{eqnarray}
for $i \in I$, $\mu \in Q$. 
For any $\l \in \bbZ_l^n$, 
we regard $\Vfin$ as a $(\Ufin)^{'}$-module through these 
homomorphisms.
Then the evaluation $\tUfin$-representations are defined
by the following method. 
\newtheorem{def SQEM}[def sec4]{Definition}
\begin{def SQEM} 
\label{def SQEM}
Let $\bf{a} \in \bbC^{\times}$, $\l=(\l_i)_{i \in I} \in \bbZ_l^n$. 
We set 
\begin{eqnarray}
&& \l_{\L_i}:=\sum_{j \in I}\l_j(\L_i, \L_j) \q (i \in I), 
\label{sig Ln}\\
&& {\bf{a}}^{\l}_{+}:=\bf{a}\e^{-\l_{\L_1}+\l_{\L_n}+n}, 
\q {\bf{a}}^{\l}_{-}:=\bf{a}(-1)^{n+1}\e^{\l_{\L_1}-\l_{\L_n}+2n+1}. 
\label{sig cl}
\end{eqnarray}
We regard $\Vfin$ as a $\tUfin$-module by using 
$(\rm{ev}_{{\bf{a}}^{\l}_{\pm}}^{\rm{fin}})^{\pm}$ 
and denote it by $\Vfin_{\bf{a}}^{\pm}$.
\end{def SQEM}
Since $\Vfin$ is irreducible as a $\Ufin$-module, 
$\Vfin_{\bf{a}}^{\pm}$ 
is a finite dimensional irreducible $\tUfin$-module of type \bf{1}. 
Thus, by Theorem \ref{thm ARSCT}, 
there exists a unique 
$\bf{P}_{\bf{a}}^{\pm}=(P_{i,\bf{a}}^{\pm})_{i \in I} \in \bbC_l[t]^n$ 
such that $\Vfin_{\bf{a}}^{\pm} \cong 
\tV_{\e}^{\rm{fin}}(\bf{P}_{\bf{a}}^{\pm})$ as a $\tUfin$-module. 
Let $i \in I$ such that $P_{i,\bf{a}}^{\pm} \neq 0$. 
Then, by Lemma \ref{lem ER}, 
there exist $\bf{a}_{(\pm,i)} \in \bbC^{\times}$ 
such that $P_{i,\bf{a}}^{(\pm)}(t)
=\prod_{p=1}^{\l_i}(t-\bf{a}_{(\pm,i)}\e^{\l_i+1-2p})$. 
Around $t=0$, we have 
\begin{eqnarray*}
 \e^{\l_i}+(\e^{\l_i}-\e^{-\l_i})
\sum_{m=1}^{\infty}(\bf{a}_{(\pm,i)}^{-1}\e^{\l_i-1}t)^m
&=&\frac{\e^{\l_i}-\bf{a}_{(\pm,i)}^{-1}\e^{-1}t}
{1-\bf{a}_{(\pm,i)}^{-1}\e^{\l_i-1}t}
=\e^{\l_i}\frac{\e^{-2\l_i}
(t-\bf{a}_{(\pm,i)}\e^{\l_i+1})}{t-\bf{a}_{(\pm,i)}\e^{-\l_i+1}}\\
&=&\e^{\rm{deg}(P^{\pm}_{i,\bf{a}})}
\frac{P^{\pm}_{i,\bf{a}}(\e^{-2}t)}{P^{\pm}_{i,\bf{a}}(t)},
\end{eqnarray*}
(see \cite{CP91} Corollary 4.2). 
Thus, by Definition \ref{def fin dim HW} (d), we obtain 
\begin{eqnarray}
\sum_{m=0}^{\infty}\wt{\L}^{\rm{fin}}_{\bf{P}^{\pm}_{\bf{a}}}
(\psi_{i,m}^{+})t^m
=\e^{\rm{deg}(P^{\pm}_{i,\bf{a}})}
\frac{P^{\pm}_{i,\bf{a}}(\e^{-2}t)}{P^{\pm}_{i,\bf{a}}(t)}
=\e^{\l_i}+(\e^{\l_i}-\e^{-\l_i})
\sum_{m=1}^{\infty}(\bf{a}_{(\pm,i)}^{-1}\e^{\l_i-1}t)^m.
\label{fac DP2}
\end{eqnarray} 
 Hence we can calculate $\bf{a}_{(\pm,i)}$ explicitly 
by the computation of 
$\wt{\L}^{\rm{fin}}_{\bf{P}^{\pm}_{\bf{a}}}(\psi_{i,1}^{+})$.
Therefore, by the similar way to the proof of 
\cite{CP94a} Theorem 3.5, we obtain the following theorem. 
\newtheorem{thm DP}[def sec4]{Theorem}
\begin{thm DP} 
\label{thm DP}
For $\l=(\l_i)_{i \in I} \in \bbZ_l^n, \bf{a} \in \bbC^{\times}$, 
let $\bf{P}^{\pm}_{\bf{a}}=(P^{\pm}_{i,\bf{a}})_{i \in I} \in \bbC_l[t]^n$ 
such that $\tV_{\e}^{\rm{fin}}(\bf{P}^{\pm}_{\bf{a}})
 \cong \Vfin_{\bf{a}}^{\pm}$. 
Then, for any $i \in I$ such that $P^{\pm}_{i,\bf{a}} \neq 1$,
we obtain
\begin{eqnarray}
 P_{i,\bf{a}}^{\pm}=\prod_{p=1}^{\l_i}(t-\e^{\l_i-2p+1}\bf{a}_{(\pm,i)}),
\label{def pia}
\end{eqnarray}
where 
\begin{eqnarray}
&& \bf{a}_{(\pm,i)}:=\bf{a}^{-1}\e^{\pm(\l^{(i)}+i)},
\label{sig ci}\\
&& \l^{(i)}:=\sum_{k=1}^{i-1}\l_k-\sum_{k=i+1}^n\l_k.
\label{sig li}
\end{eqnarray}
\end{thm DP}
Proof. 
We shall prove the case of $P_{i,a}^{+}$. 
We can also prove the case of $P_{i,a}^{-}$ similarly.
Let $v_{+}$ be the highest weight vector in 
$\tV_{\e}^{\rm{fin}}(\bf{P}^{+}_{\bf{a}})$.
By (\ref{pro DR}), for any $i \in I$, we have
\begin{eqnarray*}
e_0v_{+}=(-1)^{i+1}\e^{n+1}[f_{n},  \cdots [
f_{i+1},[ f_{1},  \cdots [f_{i-1},
x_{i,1}^{-}]_{\e^{-1}} \cdots ]_{\e^{-1}}k_{\t}^{-1}v_{+}.  
\end{eqnarray*}
Then we get 
\begin{eqnarray*}
e_ne_0v_{+}
&=&(-1)^{i+1}\e^{-\sum_{k \in I}\l_k+n+1}\{
e_nf_n[f_{n-1}, [ \cdots [f_{i-1},
x_{i,1}^{-}]_{\e^{-1}} \cdots ]_{\e^{-1}} \\
&& \qq -\e^{-1}[f_{n-1},[ \cdots [f_{i-1},
x_{i,1}^{-}]_{\e^{-1}} \cdots ]_{\e^{-1}}e_nf_n\}v_{+} \\
&=&(-1)^{i+1}\e^{-\sum_{k \in I}\l_k+n+1}\{
(\frac{k_{\al_n}-k_{\al_n}^{-1}}{\e-\e^{-1}})[f_{n-1},[\cdots [f_{i-1},
x_{i,1}^{-}]_{\e^{-1}} \cdots ]_{\e^{-1}} \\
&& \qq -\e^{-1}[f_{n-1}, [ \cdots [f_{i-1},
x_{i,1}^{-}]_{\e^{-1}} \cdots ]_{\e^{-1}}
(\frac{k_{\al_n}-k_{\al_n}^{-1}}{\e-\e^{-1}})\}v_{+} \\ 
&=&(-1)^{i+1}\e^{-\sum_{k \in I}\l_k+n+1}
([\l_n+1]_{\e}-\e^{-1}[\l_n]_{\e})[f_{n-1},[\cdots [f_{i-1},
x_{i,1}^{-}]_{\e^{-1}} \cdots ]_{\e^{-1}} v_{+} \\
&=&(-1)^{i+1}\e^{\l_n-\sum_{k \in I}\l_k+n+1}[f_{n-1},[\cdots [f_{i-1},
x_{i,1}^{-}]_{\e^{-1}} \cdots ]_{\e^{-1}} v_{+}.  
\end{eqnarray*}
By repeating this, we obtain 
\begin{eqnarray}
 e_i \cdots e_1e_{i+1} \cdots e_ne_0v_{+}
&=&e_i\{(-1)^{i+1}\e^{-\l_i+n+1}x_{i,1}^{-}v_{+}\} \no\\
&=&(-1)^{i+1}\e^{-\l_i+n+1}x_{i,0}^{+}x_{i,1}^{-}v_{+}.
\label{sig e01}
\end{eqnarray}
On the other hand, by (\ref{sig SQEH}), we have 
\begin{eqnarray*}
 e_0v_{+}=(\rm{ev}_{\bf{a}_{+}^{\l}}^{\rm{fin}})^{+}(e_0)v_{+}
=\bf{a}_{+}^{\l}\e^{-1}k_{\L_1}k_{\L_n}^{-1}[f_n,[ \cdots[f_2,f_1]_{\e^{-1}}
\cdots]_{\e^{-1}}v_{+}.
\end{eqnarray*} 
Since $(\L_1-\L_n, \t)=0$, we have
\begin{eqnarray*}
k_{\L_1}k_{\L_n}^{-1}[f_n,[ \cdots[f_2,f_1]_{\e^{-1}}
\cdots]_{\e^{-1}}
=[f_n,[ \cdots[f_2,f_1]_{\e^{-1}}
\cdots]_{\e^{-1}}k_{\L_1}k_{\L_n}^{-1}. 
\end{eqnarray*}
Moreover, we have 
$k_{\L_1}v_{+}=\e^{\l_{\L_1}}v_{+}$, 
$k_{\L_n}v_{+}=\e^{\l_{\L_n}}v_{+}$. 
So, by (\ref{sig cl}), we obtain
\begin{eqnarray*}
 e_0v_{+}
=\bf{a}\e^{n-1}[f_n,[ \cdots[f_2,f_1]_{\e^{-1}}\cdots]_{\e^{-1}}
v_{+}. 
\end{eqnarray*}
By the similar way to the above proof, we have 
\begin{eqnarray*}
 e_{i+1} \cdots e_ne_0v_{+}
=\bf{a}\e^{\sum_{k=i+1}\l_k+n-1}[f_i,[\cdots[f_2,f_1]_{\e^{-1}} 
\cdots]_{\e^{-1}}v_{+}.
\end{eqnarray*}
Here, we obtain 
\begin{eqnarray*}
&& e_1[f_i,[\cdots[f_2,f_1]_{\e^{-1}} 
\cdots]_{\e^{-1}}v_{+}
=[f_i,[\cdots[f_2,e_1f_1]_{\e^{-1}} 
\cdots]_{\e^{-1}}v_{+}\\
&&=\{[f_i,[\cdots[f_2,f_1e_1]_{\e^{-1}} 
\cdots]_{\e^{-1}}+[f_i,[\cdots[f_2,
\frac{k_{\al_1}-k_{\al_1}^{-1}}{\e-\e^{-1}}]_{\e^{-1}} 
\cdots]_{\e^{-1}}\}v_{+}\\
&&=\{[f_i,[\cdots[f_2,f_1]_{\e^{-1}} 
\cdots]_{\e^{-1}}e_1+[f_i,[\cdots[f_3,f_2]_{\e^{-1}} 
\cdots]_{\e^{-1}}(-\e^{-1}k_{\al_1}^{-1})\}v_{+}\\
&&=\e^{-\l_1-1}[f_i,[\cdots[f_3,f_2]_{\e^{-1}} 
\cdots]_{\e^{-1}}v_{+}.
\end{eqnarray*}
By repeating this, we get 
\begin{eqnarray}
 e_i \cdots e_1e_{i+1}\cdots e_ne_0v_{+}
&=&(-1)^{i-1}\bf{a}\e^{-\sum_{k=1}^{i-1}\l_k-(i-1)+\sum_{k=i+1}^n\l_k+n-1}
e_if_iv_{+} \no \\
&=&(-1)^{i-1}\bf{a}\e^{-\l^{(i)}-i+n}[\l_i]_{\e}v_{+}.
\label{sig e02}
\end{eqnarray}
Thus, by (\ref{sig e01}) and (\ref{sig e02}), we obtain
\begin{eqnarray*}
 \psi_{i,1}v_{+}
=(\e-\e^{-1})x_{i,0}^{+}x_{i,1}^{-}v_{+}
=\bf{a}\e^{\l_i-\l^{(i)}-i-1}(\e^{\l_i}-\e^{-\l_i})v_{+}.
\end{eqnarray*}
On the other hand, by (\ref{fac DP2}), we have
$\psi_{i,1}^{+}v_{+}
=\wt{\L}^{\rm{fin}}_{\bf{P}^{+}_{\bf{a}}}
(\psi_{i,1}^{+})v_{+}
=\bf{a}_{(+,i)}^{-1}\e^{\l_i-1}(\e^{\l_i}-\e^{-\l_i})$. 
It amount to $\bf{a}_{(+,i)}=\bf{a}^{-1}\e^{\l^{(i)}+i}$. 
\qed \\ \\ 
\q For $\l=(\l_i)_{i \in I} \in \bbZ_l^{n}$, 
we set $\rm{supp}(\l):=\{i\in I  \, |\, \l_i \neq 0\}$.
\newtheorem{pro pc=pc22}[def sec4]{Proposition}
\begin{pro pc=pc22} 
\label{pro pc=pc22}
Let $\l=(\l_i)_{i \in I} \in \bbZ_l^n$, 
$\bf{a}_{\pm} \in \bbC^{\times}$. \\
\q (a) If $\l=0$, then $\Vfin_{\bf{a}_{+}}^{+}$ is isomorphic to 
$\Vfin_{\bf{a}_{-}}^{-}$  as a $\tUfin$-module. \\
\q (b) In the case of $\l \neq 0$, 
$\Vfin_{\bf{a}_{+}}^{+}$ is isomorphic to $\Vfin_{\bf{a}_{-}}^{-}$ 
as a $\tUfin$-module 
if and only if $\bf{a}_{+}=\bf{a}_{-}\e^{2(\l^{(i)}+i)}$ 
for any $i \in \rm{supp}(\l)$. 
\end{pro pc=pc22}
Proof. 
(a) is obvious. So we shall prove (b). 
By Theorem \ref{thm ARSCT}, 
$\Vfin_{\bf{a}_{+}}^{+}$ is isomorphic to $\Vfin_{\bf{a}_{-}}^{-}$ 
if and only if $\bf{P}_{\bf{a}_{+}}^{+}=\bf{P}_{\bf{a}_{-}}^{-}$. 
By (\ref{fac DP2}) and Theorem \ref{thm DP}, 
we obtain $\bf{P}_{\bf{a}_{+}}^{+}=\bf{P}_{\bf{a}_{-}}^{-}$ if and only if 
$\bf{a}_{+}=\bf{a}_{-}\e^{2(\l^{(i)}+i)}$ for any $i \in \rm{supp}(\l)$. 
\qed 
\newtheorem{pro pc=pc222}[def sec4]{Proposition}
\begin{pro pc=pc222} 
\label{pro pc=pc222}
Let $\l=(\l_i)_{i \in I} \in \bbZ_l^n$ $(\l \neq 0)$, 
$\bf{a}_{\pm} \in \bbC^{\times}$. 
Let $i_1, \cdots, i_m \in I$
such that $\rm{supp}(\l)=\{i_1, \cdots,  i_m\}$ 
and $i_1 < \cdots <i_m$.
Then $\bf{a}_{+}=\bf{a}_{-}\e^{2(\l^{(i)}+i)}$ for any $i \in \rm{supp}(\l)$ 
if and only if the following conditions (a) and (b) hold. \\
\q (a) For any $2 \leq r \leq m$, 
\begin{eqnarray*}
 \l_{i_r} \equiv (-1)^{r-1}\l_{i_1}+(-1)^ri_1-i_r
+2\sum_{k=2}^{r-1}(-1)^{r-1+k}i_k \not\equiv 0
\q (\rm{mod } l).
\end{eqnarray*}
\q (b) 
\begin{eqnarray*}
\bf{a}_{+}=
\begin{cases}
\bf{a}_{-}\e^{2\sum_{k=1}^m(-1)^{k-1}i_k}& \rm {if $m$ is odd}, \\
\bf{a}_{-}\e^{2(\l_{i_1}+ \sum_{k=2}^m(-1)^ki_k)}&\rm {if $m$ is even}.
\end{cases}
\end{eqnarray*}
\end{pro pc=pc222}
Proof. 
We assume $\bf{a}_{+}=\bf{a}_{-}\e^{2(\l^{(i)}+i)}$ 
for any $i \in \rm{supp}(\l)$. 
Then $\e^{2(\l^{(i)}+i)}=\e^{2(\l^{(j)}+j)}$ 
for any $i,j \in \rm{supp}(\l)$.  
Hence, for $2 \leq r \leq m$, we have 
$\l^{(i_{r})}-\l^{(i_{r-1})}+i_{r}-i_{r-1} \equiv 0 
\, (\rm{mod }l)$.
By (\ref{sig li}), for $1 \leq r \leq m$, we obtain 
\begin{eqnarray*}
 \l^{(i_r)}=\sum_{k=1}^{i_r-1}\l_k-\sum_{k=i_r+1}^n\l_k
=\sum_{k=1}^{r-1}\l_{i_k}-\sum_{k=r+1}^m\l_{i_k}.
\end{eqnarray*} 
Thus, for $2 \leq r \leq m$, we get 
\begin{eqnarray*}
 \l^{(i_{r})}-\l^{(i_{r-1})}+i_{r}-i_{r-1} 
=\l_{i_{r}}+\l_{i_{r-1}}+i_{r}-i_{r-1} \equiv 0 \q (\rm{mod }l).
\label{root of unity}
\end{eqnarray*}
Hence $(-1)^r\l_{i_r}-(-1)^{r-1}\l_{i_{r-1}} \equiv 
(-1)^{r-2}i_{r-1}+(-1)^{r-1}i_r$. 
Therefore, we obtain 
\begin{eqnarray*}
(-1)^r\l_{i_r} &\equiv& -\l_{i_1}
+\sum_{k=2}^r\{(-1)^k\l_{i_k}-(-1)^{k-1}\l_{i_{k-1}}\}
\equiv -\l_{i_1}
+\sum_{k=2}^r\{(-1)^{k-1}i_k+(-1)^{k-2}i_{k-1}\}\\
&\equiv& -\l_{i_1}+i_1+(-1)^{r-1}i_r 
+2\sum_{k=2}^{r-1}(-1)^{k-1}i_k 
\q (2 \leq r \leq m).
\end{eqnarray*}
Thus we have (a). 
In particular, if $\l_{i_r}$ is as in (a), 
then $\l_{i_{r-1}}+\l_{i_{r}} \equiv i_{r-1}-i_r$ 
for any $2 \leq r \leq m$ 
and $\l^{(i)}+i \equiv\l^{(j)}+j$ for any $i,j \in \rm{supp}(\l)$. 
Hence, for $1 \leq r \leq m$, we have
$\l^{(i_r)}+i_r\equiv\l^{(1)}+i_1
\equiv -\sum_{k=2}^m\l_{i_k}+i_1$. 
If $m$ is odd, then we get 
\begin{eqnarray*}
- \sum_{k=2}^m\l_{i_k}+i_1&\equiv&
-(\l_{i_2}+\l_{i_3})- \cdots -(\l_{i_{m-1}}+\l_{i_m})+i_1\\
&\equiv&(-i_2+i_3)+\cdots +(-i_{m-1}+i_{m})+i_1
\equiv\sum_{k=1}^m(-1)^{k-1}i_k. 
\end{eqnarray*}
Similarly, we have the case that $m$ is even.
Therefore we obtain (b). 
So we can prove ``only if part'' of this proposition. 
The proof of ``if part''  follows  
the proof of ``only if part''. 
\qed
\newtheorem{rem pc=pc2}[def sec4]{Remark}
\begin{rem pc=pc2} 
\label{rem pc=pc2}
For $\l=(\l_i)_{i \in I} \in \bbZ_{+}^n$, let $V_q(\l)$ be the 
finite dimensional irreducible $U_q$-module with highest weight $\l$ 
of type $\bf{1}$. 
For $\bf{a} \in \bbC^{\times}$, 
let $V_q(\l)_{\bf{a}}^{\pm}$ 
be the evaluation representation of $V_q(\l)$ arising from 
$\rm{ev}_{\bf{a}}^{\pm}$ (see \cite{CP94a}). 
In the case that $q$ is not a root of unity, 
for any $\bf{a}_{\pm} \in \bbC^{\times}$, 
$V_q(\l)_{\bf{a}_{+}}^{+}$ is not isomorphic to 
$V_q(\l)_{\bf{a}_{-}}^{-}$ 
if $\#(\rm{supp}(\l)) >1$.
But, in the case that $q$ is a root of unity, 
there exist $\l \in \bbZ_l^n$ 
and $\bf{a}_{\pm} \in \bbC^{\times}$ such that 
$\Vfin_{\bf{a}_{+}}^{+}$ is isomorphic to $\Vfin_{\bf{a}_{-}}^{-}$ 
even if $\#(\rm{supp}(\l))>1$ by Proposition 
\ref{pro pc=pc22}, \ref{pro pc=pc222}.
\end{rem pc=pc2}
\section{Evaluation representations of non-restricted type}
\setcounter{equation}{0}
\renewcommand{\theequation}{\thesection.\arabic{equation}}
\subsection{Schnizer modules and 
evaluation representations} 
\q We fix the following notations. 
Let $N:=\frac{1}{2}n(n+1)$ be the number of the positive roots 
of $\sl_{n+1}$. 
Let $V_{N}$ be a $l^N$-dimensional $\bbC$-vector space 
and $\{v(m) \in V_{N}\, | 
\, m=(m_{i,j})_{1 \leq i \leq j \leq n} \in  \bbZ_l^N\}$ 
be a $\bbC$-basis of $V_{N}$. 
For $m \in \bbZ_l^N$, $m^{'} \in \bbZ^{N}$, 
we set $v(m+lm^{'}):=v(m)$. 
For $i,j \in I$, we define $\ep_{i,j}, \al_{i,j} \in \bbZ_l^N$ by 
\begin{eqnarray}
 \ep_{i,j}:=(\d_{i,r}\d_{j,s})_{1 \leq r \leq s \leq n} \q ( i \leq j ), 
\q \al_{i,j}:=\sum_{k=j+1}^i\ep_{k-1, n-i+k}-\sum_{k=j}^i \ep_{k,n-i+k} 
\q (j \leq i),
\label{sig alep}
\end{eqnarray}
where $\d_{i,j}$ is the Kronecker's symbol. 
For $i,j \in I$, 
$a=(a_{i,j})_{1 \leq i \leq j \leq n} \in (\bbC^{\times})^N$, 
$c=(c_{i,j})_{1 \leq i \leq j \leq n} \in \bbC^N$, 
we define 
\begin{eqnarray}
&& M_{i,j}(c):=\sum_{k=i-1}^{j-1}(c_{i,k}-c_{i-1,k})
+\sum_{k=i}^j(c_{i,k}-c_{i+1,k}) \q (i \leq j), 
\label{sig M} \\
&& N_{i,j}(c):=c_{j-1,n-i+j}-c_{j,n-i+j} \q (j \leq i), 
\label{sig N} \\
&& \mu_i(c):=\sum_{k=i-1}^n c_{i-1,k}
-2\sum_{k=i}^nc_{i,k}+\sum_{k=i+1}^n c_{i+1,k}, 
\label{sig mu} \\
&& a(c):=\prod_{1 \leq i \leq j \leq n}a_{i,j}^{c_{i,j}}, 
\label{sig ac}
\end{eqnarray} 
where $c_{i,j}=0$ if the index $(i,j)$ is out of range.
\newtheorem{thm SM}{Theorem}[section]
\begin{thm SM}[\cite{S94} Theorem 3.2, \cite{S93}]
\label{thm SM}
Let $a=(a_{i,j})_{1 \leq i \leq j \leq n} \in (\bbC^{\times})^{N}$, 
$b=(b_{i,j})_{1 \leq i \leq j \leq n} \in \bbC^{N}$, 
and $\l=(\l_i)_{i \in I} \in \bbC^n$. 
Then there exists a $\bbC$-algebra homomorphism 
$\rho:=\rho(a,b,\l): U_{\e}\arr \rm{End}(V_{N})$ such that 
for $i \in I$ and $m \in \bbZ_l^N$
\begin{eqnarray}
&& \rho(E_i)(v(m))=
\sum_{j=1}^i a(\al_{i,j})[N_{i,j}(m+b)]_{\e}v(m+\al_{i,j}),  
\label{sig Sch-E}\\
&& \rho(F_i)(v(m))=
\sum_{j=i}^n a_{i,j}[M_{i,j}(m+b)-\l_i]_{\e}v(m+\ep_{i,j}),  
\label{sig Sch-F}\\
&& \rho(K_{\al_i})(v(m))=\e^{\mu_i(m+b)+\l_i}v(m).
\label{sig Sch-K}
\end{eqnarray}
\end{thm SM}
We denote the $l^N$-dimensional $U_{\e}$-module 
associated with $(\rho(a,b,\l), V_N)$ by $V_{\e}(a,b,\l)$. 
We call $V_{\e}(a,b,\l)$ a ``Schnizer module''. 

Now, for $i \in I$ and $r=(r_1, \cdots, r_i) \in I^i$, we set 
\begin{eqnarray}
&& F_{\t_i}:=[F_i,[ \cdots [F_2, F_1]_{\e^{-1}} \cdots ]_{\e^{-1}}, 
\q E_{\t_i}:=[E_i,[ \cdots [E_2, E_1]_{\e^{-1}} \cdots ]_{\e^{-1}}, 
\label{sig EiFi}\\
&&  \ep_{r}:=\sum_{k=1}^i \ep_{k,r_k}, 
\q  \al_{r}:=\sum_{k=1}^i\al_{k,r_k}.
\label{sig epal} 
\end{eqnarray}
For $1 \leq s \leq i$, we set 
\begin{eqnarray}
&& R_{s,i}:=\{r_i^s=(r_{1,i}^s, \cdots, r_{i,i}^s) \in I^i \, | \,
r_{1,i}^s \geq \cdots \geq r_{s-1,i}^s \geq r_{s,i}^s 
<r_{s+1,i}^s < \cdots <r_{i,i}^s\}, \no \\
&& R_{s,i}^F:=\{r_i^s=(r_{1,i}^s, \cdots, r_{i,i}^s) \in R_{s,i} \, | \,
k \leq r_{k,i}^s \leq n \rm{ for $1 \leq k \leq i$}\},  \no \\
&& R_{s,i}^E:=\{r_i^s=(r_{1,i}^s, \cdots, r_{i,i}^s) \in R_{s,i} \, | \,
1 \leq r_{k,i}^s \leq k \rm{ for $1 \leq k \leq i$}\},  \no \\
&&R_i^F:=\bigsqcup_{s=1}^i R_{s,i}^F, 
\q R_i^E:=\bigsqcup_{s=1}^i R_{s,i}^E. 
\label{sig RiFRiE}
\end{eqnarray}
Moreover, for $c \in \bbC^N$, set
\begin{eqnarray}
&&C_i(c,r_i^s)
:=\sum_{k=1}^{s-1}M_{k,r_{k,i}^s}(c)
-\sum_{k=s+1}^{i}M_{k,r_{k,i}^s}(c) 
\q (r_i^s \in R_{i,s}^F), \no \\
&&D_i(c,r_i^s)
:=\sum_{k=1}^{s-1}N_{k,r_{k,i}^s}(c)
-\sum_{k=s+1}^{i}N_{k,r_{k,i}^s}(c) 
\q (r_i^s \in R_{i,s}^E). 
\label{sig CiDi} 
\end{eqnarray}
\newtheorem{lem SE}[thm SM]{Lemma}
\begin{lem SE}
\label{lem SE}
Let $a=(a_{i,j})_{1 \leq i \leq j \leq n} \in (\bbC^{\times})^N$, 
$b=(b_{i,j})_{1 \leq i \leq j \leq n} \in \bbC^N$, 
$\l=(\l_i)_{i \in I} \in \bbC^n$, 
$i \in I$, and $m=(m_{i,j})_{1 \leq i \leq j \leq n} \in \bbZ_l^N$. 
We have
\begin{eqnarray*}
&&F_{\t_i}v(m)
=\sum_{r_i^s \in R_i^F} (-1)^{i+s}a(\ep_{r_i^s})
\e^{C_i(m+b,r_i^s)-\l^{(s,i)}+1-s} 
[M_{s,r_{s,i}^s}(m+b)-\l_s]_{\e}v(m+\ep_{r_i^s}), \\
&& E_{\t_i}v(m)
=\sum_{r_i^s \in R_i^E} (-1)^{i+s}a(\al_{r_i^s})
\e^{D_i(m+b,r_i^s)+1-s}
[N_{s,1}(m+b)-\l_s]_{\e}v(m+\al_{r_i^s}), 
\end{eqnarray*}
in $V_{\e}(a,b,\l)$, where 
$\l^{(s,i)}:=\sum_{k=1}^{s-1}\l_k-\sum_{k=s+1}^i\l_k$.
\end{lem SE}
Proof. We shall prove the 
$F_{\t_i}$-case by the induction on $i$. 
We can prove the $E_{\t_i}$-case similarly. 
If $i=1$, then we have
\begin{eqnarray*}
 F_{\t_1}v(m)
=\sum_{r_1^1 \in R_1^F} a(\ep_{1,r_{1,1}^1})[
M_{1,r_{1,1}^1}(m+b)-\l_1]_{\e}v(m+\ep_{1,r_{1,1}^1}).
\end{eqnarray*}
We replace $r_{1,1}^1$ with $j$. 
Then we obtain 
\begin{eqnarray*}
 F_{\t_1}v(m)
=\sum_{j=1}^n a(\ep_{1,j})[M_{1,j}(m+b)-\l_1]_{\e}v(m+\ep_{1,j})
=F_1v(m).
\end{eqnarray*}
\q Now we assume that $i>1$ and we get the case of $(i-1)$. 
For $i\leq j \leq n$, $r_{i-1}^s \in R_{i-1}^F$, we set 
\begin{eqnarray*}
&& \bf{M}(r_{i-1}^s,j):=[M_{s,r_{s,i-1}^s}(m+b)-\l_s]_{\e}[
M_{i,j}(m+b)-\l_i+M_{i,j}(\ep_{r_{i-1}^s})]_{\e}\\
&& \qq \qq -\e^{C_{i-1}(\ep_{i,j},r_{i-1}^s)-1}[
M_{s,r_{s,i-1}^s}(m+b)-\l_s+M_{s,r_{s,i-1}^s}(\ep_{i,j})]_{\e}[
M_{i,j}(m+b)-\l_i]_{\e}.
\end{eqnarray*}
Then, by the assumption of the induction, we have 
\begin{eqnarray*}
&&\hspace{-30pt} F_{\t_i}v(m)
=[F_i, F_{\t_{i-1}}]_{\e^{-1}}v(m)\\
&&\hspace{-30pt}=\sum_{j=i}^n \sum_{r_{i-1}^s \in R_{i-1}^F}
(-1)^{i+s-1}a_{i,j}a(\ep_{r_{i-1}^s}) 
\e^{C_{i-1}(m+b,r_{i-1}^s)-\l^{(s,i-1)}+1-s}
\bf{M}(r_{i-1}^s,j)v(m+\ep_{r_{i-1}^s}+\ep_{i,j}).
\end{eqnarray*}
\q Now we set 
\begin{eqnarray*}
\xi(j>j^{'}):= 
\begin{cases}
1&\rm{if $j > j^{'}$}\\
0&\rm{if $j \leq j^{'}$},
\end{cases} 
\q \xi(j \leq j^{'}):= 
\begin{cases}
1&\rm{if $j \leq  j^{'}$}\\
0&\rm{if $j < j^{'}$}.
\end{cases} 
\end{eqnarray*}
Then, for any $1 \leq i \leq j \leq n$, 
$1 \leq i^{'} \leq j^{'} \leq n$, we get 
\begin{eqnarray*}
 M_{i,j}(\ep_{i^{'},j^{'}})
=-\d_{i-1,i^{'}}\xi(j>j^{'})+\d_{i,i^{'}}\xi(j>j^{'})
+\d_{i,i^{'}}\xi(j \geq j^{'})-\d_{i+1,i^{'}}\xi(j \geq j^{'}),
\end{eqnarray*}
Hence, for any $ i\leq j \leq n$, $1 \leq s \leq i-1$, 
$r_{i-1}^s \in R_{i-1}^F$, we have
\begin{eqnarray*}
 &&M_{i,j}(\ep_{r_{i-1}^s})=-\xi(j>r_{i-1,i-1}^s), 
\q M_{s,r_{s,i-1}^s}(\ep_{i,j})=-\d_{s,i-1}\xi(r_{i-1,i-1}^{i-1} \geq j), \\
&& C_{i-1}(\ep_{i,j}, r_{i-1}^s)=\xi(i-1>s)\xi(r_{i-1,i-1}^s \geq j). 
\end{eqnarray*}
Thus, we have   
\begin{eqnarray*}
 &&\bf{M}(r_{i-1}^s,j)
=[M_{s,r_{s,i-1}^s}(m+b)-\l_s]_{\e}[
M_{i,j}(m+b)-\l_i-\xi(j>r_{i-1,i-1}^s)]_{\e}\\
&& -\e^{\xi(i-1>s)\xi(r_{i-1,i-1}^s \geq j)-1}[
M_{s,r_{s,i-1}^s}(m+b)-\l_s-\d_{s,i-1}\xi(r_{i-1,i-1}^{i-1} \geq j)]_{\e}
[M_{i,j}(m+b)-\l_i]_{\e}. 
\end{eqnarray*}
Since
$[c]_{\e}-\e^{-1}[c-1]_{\e}=\e^{c-1}$ 
and $[c-1]_{\e}-\e^{-1}[c]_{\e}=-\e^{-c}$
for any $c \in \bbC$, 
we have
\begin{eqnarray*}
 &&\bf{M}(r_{i-1}^s,j)=
-\e^{M_{i,j}(m+b)+\l_i}[M_{s,r_{s,i-1}^s}(m+b)-\l_s]_{\e} 
\q (s \leq i-1, \, r_{i-1,i-1}^s<j), \\
&&\bf{M}(r_{i-1}^s,j)=0 \q (s<i-1, \, j \leq r_{i-1,i-1}^s), \\
&&\bf{M}(r_{i-1}^{i-1},j)=\e^{\wt{M}(r_{i-1}^{i-1},j)}
[M_{i,j}(m+b)-\l_i]_{\e}, 
 \q (j \leq r_{i-1,i-1}^{i-1}),
\end{eqnarray*}
where 
$\wt{M}(r_{i-1}^{i-1},j)
:=M_{i-1,r_{i-1,i-1}^{i-1}}(m+b)-\l_{i-1}-1$.
Therefore we obtain 
\begin{eqnarray*}
F_{\t_i}v(m)
&=&\sum_{j>r_{i-1,i-1}^s}\sum_{r_{i-1}^s \in R_{i-1}^F}
(-1)^{i+s}a(\ep_{r_{i-1}^s}+\ep_{i,j})
\e^{C_{i-1}(m+b,r_{i-1}^s)-M_{i,j}(m+b)-\l^{(s,i)}+1-s} \\
&& \qq \qq \qq \qq \qq \qq \qq  
[M_{s,r_{s,i-1}^s}(m+b)-\l_s]_{\e}v(m+\ep_{r_{i-1}^s}+\ep_{i,j})\\
&+&\sum_{j \leq r_{i-1,i-1}^s}\sum_{r_{i-1}^{i-1} \in R_{i-1,i-1}^F}
(-1)^{i+i}a(\ep_{r_{i-1}^{i-1}}+\ep_{i,j})
\e^{C_{i-1}(m+b,r_{i-1}^{i-1})-M_{i,j}(m+b)-\l^{(i,i)}+1-i}\\
&& \qq \qq \qq \qq \qq \qq \qq  
 [M_{i,j}(m+b)-\l_i]_{\e}v(m+\ep_{r_{i-1}^{i-1}}+\ep_{i,j}).
\end{eqnarray*}
Here, if we set 
\begin{eqnarray*}
r_i^s=(r_{1,i}^s, \cdots, r_{i,i}^s):= 
\begin{cases} 
(r_{1,i-1}^s, \cdots, r_{i-1,i-1}^s,j)&
\rm{if $s \leq i-1$ and $j>r_{i-1,i-1}^s$} \\
(r_{1,i-1}^{i-1}, \cdots, r_{i-1,i-1}^{i-1},j)& 
\rm{if $s=i$ and $j \leq r_{i-1,i-1}^{i-1}$}, 
\end{cases}
\end{eqnarray*}
then we have $F_{\t_i}$-case. 
\qed 

For $s \in I$, we set 
\begin{eqnarray}
&& R_{s}:=\{r^s=(r_k^s)_{k \in I} \in I^n \, | \,
r_{1}^s \geq \cdots \geq r_{s-1}^s \geq r_{s}^s 
<r_{s+1}^s < \cdots <r_{n}^s\}, \no \\
&& R_s^F:=\{r^s=(r_i^s)_{i \in I} \in R_s \, | \, 
k \leq r_k^s \leq n \rm{ for any $k \in I$}\}, \no\\
&& R_s^E:=\{r^s=(r_i^s)_{i \in I} \in R_s \, | \, 
1 \leq r_k^s \leq k \rm{ for any $k \in I$}\}, \no\\
&&R^F:=\bigsqcup_{s=1}^nR_s^F, 
\q R^E:=\bigsqcup_{s=1}^nR_s^E.
\label{sig RFRE}
\end{eqnarray}
Note if $r^s=(r_k^s)_{k \in I} \in R_s^F$ (resp. $R_s^E$), 
then $r_k^s=k$ for any $s \leq k \leq n$ 
(resp. $r_k^s=1$ for any $1 \leq k \leq s$).
Moreover, for $c \in \bbC^N$, we set 
\begin{eqnarray}
&&C(c,r^s):=c_{s-1,s-1}-c_{n,n}+\sum_{k=1}^nc_{1,k}
+\sum_{k=1}^{s-1}\sum_{p=r_{k+1}^s}^{r_k^s-1}c_{k,p} 
 -\sum_{k=1}^{s}\sum_{p=r_{k}^s+1}^{r_{k-1}^s}c_{k,p} 
\q (r^s \in R_s^F), \no\\
&&D(c,r^s):=-\sum_{k=n-s+2}^n c_{1,k}
-\sum_{k=s+1}^{n}(c_{r_k^s-1,n-k+r_k^s}-c_{r_k^s,n-k+r_k^s}) 
\q (r^s \in R_s^E),
\label{sig CD}
\end{eqnarray}
where $r_0^s:=n$. 
Then, by Lemma \ref{lem SE}, 
we obtain the following lemma.
\newtheorem{lem SE2}[thm SM]{Lemma}
\begin{lem SE2}
\label{lem SE2}
Let $a=(a_{i,j})_{1 \leq i \leq j \leq n} \in (\bbC^{\times})^N$, 
$b=(b_{i,j})_{1 \leq i \leq j \leq n} \in \bbC^N$, 
$\l=(\l_i)_{i \in I} \in \bbC^n$, and 
$m=(m_{i,j})_{1 \leq i \leq j \leq n} \in \bbZ_l^N$. 
We have
\begin{eqnarray*}
&& F_{\t_n}v(m)
=\sum_{r^s \in R^F}(-1)^{s+n}a(\ep_{r^s})
\e^{C(m+b,r^s)-\l^{(s)}+1-s} \\
&& \qq \qq \qq \qq \qq
[-m_{s-1,s-1}+m_{s,s}-b_{s-1,s-1}+b_{s,s}-\l_s]_{\e}
v(m+\ep_{r^s}), \\
&&E_{\t_n}v(m)
=\sum_{r^s \in R^E}(-1)^{s+n}a(\al_{r^s})
\e^{D(m+b,r^s)+1-s}[-m_{1,n-s+1}-b_{1,n-s+1}]_{\e}v(m+\al_{r^s}), 
\end{eqnarray*}
in $V_{\e}(a,b,\l)$, 
where $\l^{(s)}$ as in (\ref{sig li}).
\end{lem SE2}
Let $U_{\AA}^{'}$ (resp. $U_{\e}^{'}$) be the extended algebra 
of $U_{\AA}$ (resp. $U_{\e}$) defined by 
replacing $\{K_{\mu} \, | \mu \in Q\}$ with $\{K_{\mu} \, | \mu \in P\}$
(see Definition \ref{def QA}). 
By (\ref{pro EH2}), we have 
$\rm{ev}_{\bf{a}}^{\pm}(\tU_{\AA}) \subset U_{\AA}^{'}$ 
($\bf{a} \in \bbC^{\times}$).
So we obtain the evaluation homomorphisms 
$\rm{ev}_{\bf{a}}^{\pm}: \tue \arr U_{\e}^{'}$ 
as in Proposition \ref{pro EH}. \\
\q On the other hand, 
by (\ref{def fundamental weights}), 
we can regard an arbitrary Schnizer module $V_{\e}(a, b, \l)$ 
as a $U_{\e}^{'}$-module if we define  
\begin{eqnarray}
K_{\L_i}v(m)
:=\e^{-\sum_{k=i}^n(m_{i,k}+b_{i,k})+\l_{\L_i}}v(m),
\label{def KLi}
\end{eqnarray}
for any $i \in I$, $m \in \bbZ_l^N$, 
where $\l_{\L_i}$ as in (\ref{sig Ln}).
\newtheorem{def SE}[thm SM]{Definition}
\begin{def SE}
\label{def SE}
Let $a=(a_{i,j})_{1 \leq i \leq j \leq n} \in (\bbC^{\times})^{N}$, 
$b=(b_{i,j})_{1 \leq i \leq j \leq n} \in \bbC^{N}$, 
$\l=(\l_i)_{i \in I} \in \bbC^n$, 
and $\bf{a} \in \bbC^{\times}$. 
Then we define $\wt{\rm{ev}}_{\bf{a}}^{\pm}
:=\wt{\rm{ev}}_{\bf{a}}^{\pm}(a,b,\l)
:=\rho(a,b,\l) \circ \rm{ev}^{\pm}_{\bf{a}^{\l}_{\pm}}
: \tUe \arr \rm{End}(V_N)$, 
where $\rho=\rho(a,b,\l)$ is as in Theorem \ref{thm SM} 
and $\bf{a}^{\l}_{\pm}$ are as in (\ref{sig cl}). 
We denote the $l^N$-dimensional $\tU_{\e}$-module 
associated with $(\wt{\rm{ev}}_{\bf{a}}^{\pm}, V_N)$ 
by $V_{\e}(a,b,\l)_{\bf{a}}^{\pm}$. 
\end{def SE}
For $c \in C^N$, we set 
\begin{eqnarray}
&&C_{E}(c,r^s)
:=c_{s-1,s-1}+\sum_{k=1}^{s-1}\sum_{p=r_{k+1}^s}^{r_k^s-1}c_{k,p}
  -\sum_{k=1}^s\sum_{p=r_k^s+1}^{r_{k-1}^s}c_{k,p}  
\q (r^s \in R_s^F),\no\\
&&D_{F}(c,r^s)
:=-c_{n,n}+\sum_{k=1}^{n-s+1}c_{1,k}
-\sum_{k=s+1}^{n}
(c_{r_{k}^s-1,n-k+r_k^s}-c_{r_k^s,n-k+r_k^s}) 
\q (r^s \in R_s^E). \no \\
\label{sig CEDF}
\end{eqnarray}
\newtheorem{pro SE}[thm SM]{Theorem}
\begin{pro SE}
\label{pro SE}
Let $a=(a_{i,j})_{1 \leq i\leq j \leq n} \in (\bbC^{\times})^N$, 
$b=(b_{i,j})_{1 \leq i\leq j \leq n} \in \bbC^N$, 
$\l=(\l_i)_{i \in I} \in \bbC^n$, 
and $\bf{a} \in \bbC^{\times}$. 
Then, for any $i \in I$ and $m \in \bbZ_l^N$, we obtain 
\begin{eqnarray*}
&&\wt{\rm{ev}}^{\pm}_{\bf{a}}(E_i)(v(m))=\rho(E_i)(v(m)), 
\q \wt{\rm{ev}}^{\pm}_{\bf{a}}(F_i)(v(m))=\rho(F_i)(v(m)), \\
&& \wt{\rm{ev}}^{\pm}_{\bf{a}}(K_{\al_i})(v(m))=\rho(K_{\al_i})(v(m)),  \\
&&\wt{\rm{ev}}^{\pm}_{\bf{a}}(E_0)(v(m))
=\bf{a}\sum_{r^s \in R^F}(-1)^{s+n}a(\ep_{r^s})
\e^{\pm(C_{E}(m+b,r^s)-\l^{(s)}-s)+n}\\
&&\qq \qq \qq \qq \qq \qq [
-m_{s-1,s-1}+m_{s,s}-b_{s-1,s-1}+b_{s,s}-\l_s]_{\e}v(m+\ep_{r^s}), \\
&&\wt{\rm{ev}}^{\pm}_{\bf{a}}(F_0)(v(m))
={\bf{a}}^{-1}\sum_{r^s \in R^E}(-1)^{s-1}a(\al_{r^s})
\e^{\pm(D_{F}(m+b,r^s)-s+n+1)-n} \\
&&\qq \qq \qq \qq \qq \qq
[-m_{1,n-s+1}-b_{1,n-s}]_{\e}v(m+\al_{r^s}),
\end{eqnarray*}
where $\ep_{r^s}, \al_{r^s}$ as in (\ref{sig epal}), 
$R^F, R^E$ as in (\ref{sig RFRE}), 
and $C_E(m+b, r^s), D_F(m+b, r^s)$ as in (\ref{sig CEDF})
\end{pro SE}
Proof. 
By Proposition \ref{pro EH}, Theorem \ref{thm SM}, 
Lemma \ref{lem SE2}, and (\ref{def KLi}), 
we obtain the $\wt{\rm{ev}}_{\bf{a}}^{+}$-case. 
Similarly, we obtain the $\wt{\rm{ev}}_{\bf{a}}^{-}$-case.
\qed 
\subsection{Nilpotent modules}
\newtheorem{def NM}[thm SM]{Definition}
\begin{def NM} 
\label{def NM}
Let $V$ be a $\tU_{\e}$-module (resp. $U_{\e}$-module). 
We assume $E_{\b}^l=F_{\b}^l=E_{(i,sl\d)}=F_{(i,sl\d)}=0$ on $V$ 
for any $\b \in \tD_{+}^{\rm{re}}$, $i \in I$, $s \in \bbN$ 
(resp. $\bar{E}_{\ga}^l=\bar{F}_{\ga}^l=0$ on $V$ 
for any $\ga \in \D_{+}$). 
Then we call $V$  
``nilpotent'' $\tU_{\e}$-module (resp. $U_{\e}$-module). 
In particular, if $K_{\mu}^l=1$ on $V$ for any $\mu \in Q$, 
then we call $V$  
nilpotent $\tU_{\e}$-module (resp. $U_{\e}$-module) of ``type \bf{1}''.
\end{def NM}
For $\l \in \bbZ_l^n$, 
let $\Vfin$ be the $\Ufin$-module in \S 4.1. 
By Proposition \ref{pro ISO}, 
we can regard $\Vfin$ as a nilpotent $U_{\e}$-module. 
We denote $U_{\e}$-module $\Vfin$ by $\Vnil$. 
By Proposition \ref{pro ISO} and Proposition \ref{thm RSCT}, 
we obtain the following proposition.
\newtheorem{pro NCT}[thm SM]{Proposirion}
\begin{pro NCT}
\label{pro NCT}
For any $\l \in \bbZ_l^n$, 
$\Vnil$ is a finite dimensional irreducible nilpotent $U_{\e}$-module 
of type $\bf{1}$. 
Conversely, for any finite dimensional irreducible nilpotent 
$U_{\e}$-module $V$ of type \bf{1}, 
there exists a unique $\l \in \bbZ_l^n$ such that 
$V$ is isomorphic to $\Vnil$. 
\end{pro NCT}
We can construct $\Vnil$ as a $U_{\e}$-submodule of 
Schnizer module $V_{\e}(a,b,\l)$ as follows 
(\cite{AN},\cite{N}). 
For $i,j \in I$ $(i \leq j)$, $\l \in \bbZ_l^n$, we set 
\begin{eqnarray}
&&a_{i,j}^{(0)}:=1, \q b_{i,j}^{(0)}:=0, 
\q a^{(0)}:=(a_{i,j}^{(0)})_{1 \leq i \leq j \leq n}, 
\q b^{(0)}:=(b_{i,j}^{(0)})_{1 \leq i \leq j \leq n}, 
\label{sig a0b0} \\
&&\rho_{\l}^{0}:=\rho(a^{(0)}, b^{(0)}, \l), 
\q V_{\e}^0(\l):=V_{\e}(a^{(0)}, b^{(0)}, \l). 
\label{sig rho0V0}
\end{eqnarray}
We denote $v(0)$ in $V_{\e}^0(\l)$ by $v_{\l}(0)$. 
For $\l=(\l_i)_{i \in I} \in \bbZ_l^n$, 
we define $m^{\l}=(m_{i,j}^{\l})_{1 \leq i \leq j \leq n} 
\in \bbZ_l^n$ by 
\begin{eqnarray}
 m_{i,j}^{\l} \equiv \sum_{k=1}^i\l_{j-k+1} \q \rm{(mod $l$)}
\q 1 \leq i \leq j \leq n.
\label{sig lowest}
\end{eqnarray}
\newtheorem{pro SMPV}[thm SM]{Proposition}
\begin{pro SMPV}
\label{pro SMPV}
Let $\l \in \bbZ_l^n$ and $v \in V_{\e}^0(\l)$. \\
\q (a) We have $E_iv=0$ for any $i \in I$ 
if and only if $v \in \bbC v_{\l}(0)$. \\
\q (b) We have $F_iv=0$ for any $i \in I$ 
if and only if $v \in \bbC v(m^{\l})$. 
\end{pro SMPV}
Proof. 
By \cite{AN}, we obtain (a). So we shall prove (b). \\
\q ``If part''. 
By (\ref{sig lowest}), we have 
\begin{eqnarray}
m^{\l}_{i,i} -m^{\l}_{i-1,i-1}=\l_i,
\q m^{\l}_{i,j}-m^{\l}_{i-1,j}=\l_{j-i+1},
\q m^{\l}_{i,j}-m^{\l}_{i+1,j}=-\l_{j-i},
\label{mlnoseisitu}
\end{eqnarray}
for any $1 \leq i \leq j \leq n$.
Hence, by (\ref{sig M}), we get
\begin{eqnarray*}
 M_{i,j}(m^{\l})
&=&m_{i,i}^{\l}-m^{\l}_{i-1,i-1}
+\sum_{k=i}^{j-1}(m^{\l}_{i,k}-m^{\l}_{i-1,k})
+\sum_{k=i+1}^j(m^{\l}_{i,k}-m^{\l}_{i+1,k})\\
&=&\l_i+(\sum_{k=i}^{j-1}\l_{k-i+1}-\sum_{k=i+1}^j\l_{k-i})=\l_i,\\
\end{eqnarray*}
for any $1 \leq i \leq j \leq n$.
Therefore, by (\ref{sig Sch-F}), we obtain 
\begin{eqnarray*}
 F_iv(m^{\l})
=\sum_{j=i}^n[M_{i,j}(m^{\l})-\l_i]_{\e}v(m^{\l}+\ep_{i,j})
=\sum_{j=i}^n[\l_i-\l_i]_{\e}v(m^{\l}+\ep_{i,j})=0 
\q (i \in I).
\end{eqnarray*}
\q ``Only if part''. 
Let $v=\sum_{m \in \bbZ_l^n}c_mv(m) \in V(\l) (c_m \in \bbC)$. 
We assume that $F_iv=0$ for any $i \in I$.
Set 
\begin{eqnarray*}
 \bbZ_l^n(r):=\{m=(m_{i,j})_{1 \leq i\leq j \leq n} \in \bbZ_l^n \, | \, 
m_{i,j}=m_{i,j}^{\l} \rm{ if } j-i<r\}\,\,
(r=1,\cdots,n),\,\,
\bbZ_l^n(0):=\bbZ_l^n.
\end{eqnarray*} 
Then we have 
\begin{eqnarray*}
 \bbZ_l^n=\bbZ_l^n(0) \supset \bbZ_l^n(1) \supset \cdots \supset 
\bbZ_l^n(n)=\{m^{\l}\}.
\end{eqnarray*} 
We shall prove that there exist $c_m$ ($m\in \bbZ_l^n(r)$)
such that $v=\sum_{m \in \bbZ_l^n(r)}c_mv(m)$ 
by the induction on $r$. 
Indeed, if we can prove this claim, then we obtain 
$v=c_{m^{\l}}v(m^{\l}) \in \bbC v(m^{\l})$. 
If $r=0$, then there is nothing to prove. 
So we assume that $r >0$ and 
we obtain the case of $(r-1)$. \\
\q By the similar way to the proof of ``if part'', 
for any $m \in \bbZ_l^n(r-1)$,
we have $M_{i,j}(m)=\l_i$ if $j-i<r-1$. 
Moreover, for any $m \in \bbZ_l^n(r-1)$ and $i \in I$ 
such that $i+r-1 \leq n$, we get 
\begin{eqnarray*}
 M_{i,i+r-1}(m)
&=&\sum_{k=i}^{i+r-2}m_{i,k}
-\sum_{k=i-1}^{i+r-2}m_{i-1,k}
+\sum_{k=i}^{i+r-1}m_{i,k}
-\sum_{k=i+1}^{i+r-1}m_{i+1,k}\\
&=&M_{i,i+r-1}(m^{\l})-m_{i-1,i+r-2}+m_{i,i+r-1}
+m_{i-1,i+r-2}^{\l}-m_{i,i+r-1}^{\l}\\
&=&\l_i-m_{i-1,i+r-2}+m_{i,i+r-1}-\l_{i+r-1}.
\end{eqnarray*}
Therefore, by (\ref{sig Sch-F}), we get 
\begin{eqnarray}
 &&F_{i}v=\sum_{m \in \bbZ_l^n(r-1)}
c_m[m_{i,i+r-1}-m_{i-1,i+r-2}-\l_{i+r-1}]_{\e}v(m+\ep_{i,i+r-1}) \no\\
&&\qq  +\sum_{j=i+r}^{n}\sum_{m \in \bbZ_l^n(r-1)}
c_m[M_{i,j}(m)-\l_i]_{\e}v(m+\ep_{i,j})=0
 \q (i\leq n-r+1). 
\label{001}
\end{eqnarray}
\q Now, for $1 \leq s \leq n-r+1$, we set 
\begin{eqnarray*}
 \bbZ_l^n(r-1,s)&:=&\{m=(m_{i,j})_{1 \leq i \leq j \leq n} 
\in \bbZ_l^n(r-1)| \\
&&\qq  m_{i,i+r-1}-m_{i-1,i+r-2} \equiv\l_{i+r-1} 
 \rm{ (mod $l$)} \rm{ if $s \leq i \leq n-r+1$} \}.
\end{eqnarray*} 
We have $m_{i,i+r-1}-m_{i-1,i+r-2} \equiv \l_{i+r-1}$ 
for any $1 \leq i \leq n-r+1$ 
if and only if $m_{i,i+r-1}=m_{i,i+r-1}^{\l}$ 
for any $1 \leq i \leq n-r+1$.
Hence we get 
\begin{eqnarray*}
 \bbZ_l^n(r-1) \supset \bbZ_l^n(r-1,n-r+1) \supset \cdots \supset 
\bbZ_l^n(r-1,1)=\bbZ_l^n(r).
\end{eqnarray*}
So it is enough to prove that there exist 
$c_m$ ($m\in \bbZ_l^n(r-1,s)$)
such that
$v=\sum_{m \in \bbZ_l^n(r-1,s)}c_mv(m)$ 
for any $1 \leq s \leq n-r+1$. 
We shall prove this claim by the induction on $s$.
By (\ref{001}), we obtain 
\begin{eqnarray*}
 &&F_{n-r+1}v=\sum_{m \in \bbZ_l^n(r-1)}
c_m[m_{n-r+1,n}-m_{n-r,n-1}-\l_{n}]_{\e}v(m+\ep_{n-r+1,n})=0.
\end{eqnarray*}
Thus, $c_m[m_{n-r+1,n}-m_{n-r,n-1}-\l_{n}]_{\e}=0$ 
for any $m \in \bbZ_l^n(r-1)$. 
So if $c_m \neq 0$ for any $m \in \bbZ_l^n(r-1)$, 
then $m_{n-r+1,n}-m_{n-r,n-1} \equiv\l_{n}$.
Hence $v=\sum_{m \in \bbZ_l^n(r-1,n-r+1)}c_mv(m)$. \\
\q Now we assume that $s<n-r+1$ and we obtain the case of $(s+1)$. 
By (\ref{001}), we get 
\begin{eqnarray*}
 F_sv&=&\sum_{m \in \bbZ_l^n(r-1,s+1)}
c_m[m_{s,s+r-1}-m_{s-1,s+r-2}-\l_{s+r-1}]_{\e}v(m+\ep_{s,s+r-1}) \no\\
&&\qq \qq +\sum_{j=s+r}^{n}\sum_{m \in \bbZ_l^n(r-1,s+1)}
c_m[M_{s,j}(m)-\l_s]_{\e}v(m+\ep_{s,j})=0.
\end{eqnarray*}
Here, for $m=(m_{i,j})_{1 \leq i \leq j \leq n} 
\in \bbZ_l^n(r-1,s+1)$, we obtain 
\begin{eqnarray*}
 (m+\ep_{s,s+r-1})_{s+1,s+r}-(m+\ep_{s,s+r-1})_{s,s+r-1}
 &\equiv&(m_{s+1,s+r}-m_{s,s+r-1})-1 \\
&\equiv&\l_{s+r}-1,\\
 (m+\ep_{s,j})_{s+1,s+r}-(m+\ep_{s,j})_{s,s+r-1}
&\equiv&m_{s+1,s+r}-m_{s,s+r-1}\\
&\equiv&\l_{s+r} \q (s+r \leq j \leq n).
\end{eqnarray*}
Hence, by the linearly independence, we obtain 
$c_m[m_{s,s+r-1}-m_{s-1,s+r-2}-\l_{s+r-1}]_{\e}=0$ for any $m \in
\bbZ_l^n(r-1, s+1)$. 
So if $c_m \neq 0$, 
then $m_{s,s+r-1}-m_{s-1,s+r-2} \equiv\l_{s+r-1}$ 
for any $m \in \bbZ_l^n(r-1,s+1)$. 
Then we have $v=\sum_{m \in \bbZ_l^n(r-1,s)}c_mv(m)$. 
\qed \\ \\
\q For $\l \in \bbZ_l^n$, 
let $V_{\e}^0(\l)$ be as in (\ref{sig rho0V0}) and
$\Lnil$ the $U_{\e}$-submodule of 
$V_{\e}^0(\l)$ generated by $v_{\l}(0)$. 
\newtheorem{pro NM}[thm SM]{Theorem}
\begin{pro NM}[\cite{AN},\cite{N}] 
\label{pro NM}
For any $\l \in \bbZ_l^n$, 
$L_{\e}^{\rm{nil}}(\l)$ is isomorphic to 
$\Vnil$ as a $U_{\e}$-module.
\end{pro NM}
For $\bf{a} \in \bbC^{\times}$, $\l \in \bbZ_l^n$, 
let $V_{\e}^0(\l)_{\bf{a}}^{\pm}$ 
(resp. $\Lnil_{\bf{a}}^{\pm}$) 
be the evaluation representation of 
$V_{\e}^0(\l)$ (resp. $\Lnil_{\bf{a}}$) 
(see Definition \ref{def SE}). 
Then $\Lnil_{\bf{a}}^{\pm}$ is the $\tUe$-submodule of 
$V_{\e}^0(\l)_{\bf{a}}^{\pm}$ generated by $v_{\l}(0)$
respectively.

Now, let $\wt{\phi}: \tU_{\e}/\tI_{\e} \arr \tUfin$ 
(resp. $\phi: U_{\e}/I_{\e} \arr \Ufin$)
be the isomorphism in Theorem \ref{thm ISO} 
(resp. Proposition \ref{pro ISO}). 
Let $I_{\e}^{'}$ be the two sided ideal of $\ue^{'}$ 
generated by $\{\bar{E}_{\ga}^l, \bar{F}_{\ga}^l, K_{\mu}^{2l}-1 
\, | \, \ga \in \D_{+}, \mu \in P\}$. 
Then we can regard $\phi$ as an isomorphism 
from $U_{\e}^{'}/I_{\e}^{'}$ to $(\Ufin)^{'}$.
Let $\wt{\pi}: \tU_{\e} \arr \tU_{\e}/\tI_{\e}$ 
(resp. $\pi: U_{\e}^{'}  \arr U_{\e}^{'}/I_{\e}^{'}$)
be the projection and 
$\rm{ev}_{\bf{a}}^{\pm}: \tUe \arr U_{\e}^{'}$ 
(resp. $(\rm{ev}_{\bf{a}}^{\rm{fin}})^{\pm}: \tUfin \arr (\Ufin)^{'}$) 
the evaluation homomorphism in 
Proposition \ref{pro EH} (resp. (\ref{sig SQEH}))). 
Then, by the definition of these maps, 
the following diagram commutes:
\begin{equation}
\begin{CD}
\tUe @> \rm{ev}_{\bf{a}}^{\pm} >> U_{\e}^{'} \\
@V \wt{\pi} VV  @VV \pi V \\
\tUe / \wt{I}_{\e}  @. U_{\e}^{'}/ I_{\e}^{'} \\
@V \wt{\phi}VV @AA \bar{\phi}^{-1}A \\
\tUfin @> (\rm{ev}^{\rm{fin}}_{\bf{a}})^{\pm} >>(\Ufin)^{'}
\end{CD}
\label{fig}
\end{equation}
\newtheorem{pro ENM}[thm SM]{Proposirion}
\begin{pro ENM} 
\label{pro ENM}
For any $\bf{a} \in \bbC^{\times}$ and $ \l \in \bbZ_l^n$, 
$\Lnil_{\bf{a}}^{\pm}$ is a finite dimensional irreducible nilpotent 
$\tU_{\e}$-module of type \bf{1}.
\end{pro ENM}
Proof. 
We shall prove the case of $\Lnil_{\bf{a}}^{+}$. 
Since we can prove the case of $\Lnil_{\bf{a}}^{-}$ similarly.
By Theorem \ref{pro NM},
$\Lnil_{\bf{a}}^{+}$ is a finite dimensional irreducible 
$\tU_{\e}$-module of type \bf{1}. 
So we shall prove that 
$\Lnil_{\bf{a}}^{+}$ is a nilpotent $\tU_{\e}$-module. \\
\q For $\l \in \bbZ_l^n$, let $\rho^0_{\l}$ as in (\ref{sig rho0V0}). 
We define $\bar{\rho}_{\l}^0: U_{\e}^{'}/I_{\e}^{'}\arr 
\rm{End}(L_{\e}^{\rm{nil}}(\l))$ by 
$\bar{\rho}_{\l}^0(u+I_{\e}):=\rho_{\l}^0(u)$ for any $u \in U_{\e}^{'}$. 
Since $\Lnil$ is a nilpotent $U_{\e}$-module, 
$\bar{\rho}_{\l}^0$ is well defined. 
Then, for any $u \in U_{\e}^{'}$, $v \in \Lnil$, we have 
$u.v=\bar{\rho}_{\l}^0 \circ \pi(u)(v)$ on $\Lnil$. 
Hence, for any $u \in \tUe$ and
 $v \in \Lnil_{\bf{a}}^{+}$, we get 
\begin{eqnarray*}
u.v=\bar{\rho}_{\l}^0 \circ \pi \circ \rm{ev}^{+}_{\bf{a}}(u)(v)
=\bar{\rho}_{\l}^0 \circ (\bar{\phi}^{'})^{-1} 
\circ (\rm{ev}_{\bf{a}}^{\rm{fin}})^{+} \circ \wt{\phi} \circ \wt{\pi}(u)(v) 
\q \rm{on $\Lnil_{\bf{a}}^{+}$},
\end{eqnarray*}
by (\ref{fig}). 
Since $\wt{\pi}(\tI_{\e})=0$, 
we obtain $\wt{I}_{\e}=0$ on $\Lnil_{\bf{a}}$.
Therefore $\Lnil_{\bf{a}}^{+}$ is a nilpotent $\tU_{\e}$-module. 
\qed \\ \\
\q By Theorem \ref{thm ISO} and Proposition \ref{pro ENM}, 
we can regard $\Lnil_{\bf{a}}^{\pm}$ as a $\tUfin$-module. 
We denote $\tUfin$-module $\Lnil_{\bf{a}}^{\pm}$ 
by $\Lfin_{\bf{a}}^{\pm}$. 
Let $\bf{P}_{\bf{a}}^{\pm}$ be as in (\ref{def pia}) 
and $\wt{V}_{\e}^{\rm{fin}}(\bf{P}^{\pm}_{\bf{a}})$ 
the evaluation representation of $\tUfin$ in \S 4.
Then, by Theorem \ref{thm DP}, \ref{pro NM}, 
$\Lfin_{\bf{a}}^{\pm}$ is isomorphic to 
$\wt{V}_{\e}^{\rm{fin}}(\bf{P}^{\pm}_{\bf{a}})$ as a $\tUfin$-module. 
Hence, by Proposition \ref{pro pc=pc22}, 
we obtain the following proposition. 
\newtheorem{pro isomorphic condition}[thm SM]{Proposition}
\begin{pro isomorphic condition} 
\label{pro isomorphic condition}
Let $\l=(\l_i)_{i \in I} \in \bbZ_l^n$,  
$\bf{a}_{\pm} \in \bbC^{\times}$. \\
\q(a) If $\l=0$, then $\Lnil_{\bf{a}_{+}}^{+}$ is isomorphic to 
$\Lnil_{\bf{a}_{-}}^{-}$ as a $\tUe$-module. \\
\q (b) In the case of $\l \neq 0$, 
$\Lnil_{\bf{a}_{+}}^{+}$ is isomorphic to $\Lnil_{\bf{a}_{-}}^{-}$ 
as a $\tUe$-module 
if and only if $\bf{a}_{+}=\bf{a}_{-}\e^{2(\l^{(i)}+i)}$ 
for any $i \in \rm{supp}(\l)$. 
\end{pro isomorphic condition}
\subsection{Alternative proof of Proposition 
\ref{pro isomorphic condition}(b)}
We can also prove Proposition \ref{pro isomorphic condition} (b)
without using the theory of 
restricted type.
We give here the alternative proof.

``Proof of only if part''. 
We assume that $\Lnil_{\bf{a}_{+}}^{+} \cong \Lnil_{{\bf{a}}_{-}}^{-}$. 
Then there exists a $\tUe$-module isomorphism 
$\phi : \Lnil_{\bf{a}_{+}}^{+} \arr \Lnil_{{\bf{a}}_{-}}^{-}$. 
By Proposition \ref{pro SMPV} (a), 
there exists $d \in \bbC^{\times}$ such that $\phi(v_{\l}(0))=dv_{\l}(0)$. 
Since $\Lnil_{\bf{a}_{+}}^{+}$ is generated by $v_{\l}(0)$ as a $\tUe$-module, 
we obtain $\phi(v)=dv$ for any $v \in \Lnil_{\bf{a}_{+}}$. 
Hence we have 
\begin{eqnarray}
\wt{\rm{ev}}^{+}_{\bf{a}_{+}}(E_0)v_{\l}(0)
=d^{-1}\phi(\wt{\rm{ev}}^{+}_{\bf{a}_{+}}(E_0)v_{\l}(0)) 
=d^{-1}\wt{\rm{ev}}^{-}_{{\bf{a}}_{-}}(E_0)\phi(v_{\l}(0))
=\wt{\rm{ev}}^{-}_{{\bf{a}}_{-}}(E_0)v_{\l}(0).
\label{e0v=e0v}
\end{eqnarray}
\q For $c=(c_{i,j})_{1 \leq i \leq j \leq n} \in \bbC^{N}$, 
$r^s \in R^{F}$, 
let $C_{E}(c,r^s)$ be as in (\ref{sig CEDF}). 
Then $C_{E}(0,r^s)=0$ for  $r^s \in R^F$. 
By Theorem \ref{pro SE}, we obtain 
\begin{eqnarray}
\wt{\rm{ev}}_{\bf{a}_{\pm}}^{\pm}(E_0)v_{\l}(0)
=\bf{a}_{\pm}\sum_{r^s \in R^F}(-1)^{s+n}\e^{\mp(\l^{(s)}+s)+n}[
-\l_s]_{\e}v(\ep_{r^s}).
\label{sig E0-actions}
\end{eqnarray}
Since $\{v(\ep_{r^s}) \,|\, r^s \in R^F\}$ is linearly independent, 
by (\ref{e0v=e0v}) and (\ref{sig E0-actions}), we have
$\bf{a}_{+}\e^{-\l^{(s)}-s}[\l_s]_{\e}
={\bf{a}}_{-}\e^{\l^{(s)}+s}[\l_s]_{\e}$ 
for any $s \in I$. 
Hence $\bf{a}_{-}={\bf{a}_{+}}\e^{2(\l^{(i)}+i)}$ 
for any $i \in \rm{supp}(\l)$. 

``Proof of if part''. 
We assume that 
${\bf{a}}_{-}=\bf{a}_{+}\e^{2(\l^{(i)}+i)}$ for any $i \in \rm{supp}(\l)$. 
By the definition of $\wt{\rm{ev}}^{\pm}_{\bf{a}_{\pm}}$, we have 
$\wt{\rm{ev}}^{+}_{\bf{a}_{+}}(E_i)
=\wt{\rm{ev}}^{-}_{\bf{a}_{-}}(E_i)$, 
$\wt{\rm{ev}}^{+}_{\bf{a}_{+}}(F_i)
=\wt{\rm{ev}}^{-}_{\bf{a}_{-}}(F_i)$, and
$\wt{\rm{ev}}^{+}_{\bf{a}_{+}}(K_{\al_i})
=\wt{\rm{ev}}^{-}_{\bf{a}_{-}}(K_{\al_i})$ 
on $\Vnil$ for any $i \in I$. 
So it is enough to prove that 
$\wt{\rm{ev}}^{+}_{\bf{a}_{+}}(E_0)=\wt{\rm{ev}}^{-}_{\bf{a}_{-}}(E_0)$ and 
$\wt{\rm{ev}}^{+}_{\bf{a}_{+}}(F_0)=\wt{\rm{ev}}^{-}_{\bf{a}_{-}}(F_0)$ 
on $\Vnil$. 
By (\ref{sig E0-actions}), 
we obtain $\wt{\rm{ev}}^{+}_{\bf{a}_{+}}(E_0)v_{\l}(0)
=\wt{\rm{ev}}^{-}_{\bf{a}_{-}}(E_0)v_{\l}(0)$. 
On the other hand, for any $j_1, \cdots j_r \in I (r \in \bbN)$, we get 
\begin{eqnarray}
&& \wt{\rm{ev}}^{+}_{\bf{a}_{+}}(E_0)(F_{j_1} \cdots F_{j_r} v_{\l}(0))
=\wt{\rm{ev}}^{+}_{\bf{a}_{+}}(F_{j_1} \cdots 
F_{j_r})(\wt{\rm{ev}}^{+}_{\bf{a}_{+}}(E_0) v_{\l}(0)) \no \\
&&\qq=\wt{\rm{ev}}^{-}_{{\bf{a}}_{-}}
(F_{j_1} \cdots F_{j_r})(\wt{\rm{ev}}^{-}_{\bf{a}_{-}}(E_0) v_{\l}(0))
 =\wt{\rm{ev}}^{-}_{\bf{a}_{-}}(E_0)(F_{j_1} \cdots F_{j_r} v_{\l}(0)). \no
\end{eqnarray}
Since $\Lnil$ is spanned by $U_{\e}^{-}v_{\l}(0)$ as a $\bbC$-vector space, 
we obtain $\wt{\rm{ev}}^{+}_{\bf{a}_{+}}(E_0)
=\wt{\rm{ev}}^{-}_{{\bf{a}}_{-}}(E_0)$. 

Now, for $c=(c_{i,j})_{1 \leq i \leq j \leq n} \in \bbC^{N}$, 
$r^s \in R^{E}$, let $D_{F}(c,r^s)$ 
be as in Proposition \ref{pro SE} 
and $m^{\l}$ be as in (\ref{sig lowest}). 
Then, for any $r^s \in R^E$, we have 
\begin{eqnarray*}
 D_{F}(m^{\l},r^s)
&=&-m_{n,n}^{\l}+\sum_{k=1}^{n-s+1}m_{1,k}^{\l}
-\sum_{k=s+1}^n(m_{r_k^s-1,n-k+r_k^s}^{\l}-m_{r_k^s,n-k+r_k^s}^{\l})\\
&=&-\sum_{k=1}^n\l_{n-k+1}+\sum_{k=1}^{n-s+1}\l_k+\sum_{k=s+1}^n\l_{n-k+1}
=\l^{(n-s+1)},
\end{eqnarray*}
(see (\ref{mlnoseisitu})).
 Hence, by Theorem \ref{pro SE}, we get 
\begin{eqnarray*}
\widetilde{\rm{ev}}^{\pm}_{\bf{a}_{\pm}}(F_0)v(m^{\l})
=\bf{a}_{\pm}^{-1}\sum_{r^s \in R^E}(-1)^{s-1}
\e^{\pm(\l^{(n-s+1)}-s+n+1)-n}[-\l_{n-s+1}]_{\e}v(m^{\l}+\al_{r^s}).
\end{eqnarray*}
By the assumption, 
if $\l_{n-s+1} \neq 0$, then we have
$\bf{a}_{+}=\bf{a}_{-}\e^{2(\l^{(n-s+1)}+(n-s+1))}$. 
So we obtain 
\begin{eqnarray*}
 \bf{a}_{+}^{-1}\e^{\l^{(n-s+1)}-s+n+1}
=\bf{a}_{-}^{-1}\e^{-2\l^{(n-s+1)}-2(n-s+1)+\l^{(n-s+1)}-s+n+1}
=\bf{a}_{-}^{-1}\e^{-(\l^{(n-s+1)}-s+n+1)}.
\end{eqnarray*}
We have that $\widetilde{\rm{ev}}^{+}_{\bf{a}^{+}}(F_0)v(m^{\l})
=\widetilde{\rm{ev}}^{-}_{\bf{a}_{-}}(F_0)v(m^{\l})$. 

On the other hand, by the similar way to the proof of 
Proposition 5.6 in \cite{N}, we obtain that 
there exists a nonzero vector $v_{L}\in \Lnil$ 
such that $F_iv_{L}=0$ for any $i \in I$. 
Hence, by Proposition \ref{pro SMPV} (b), we obtain $v(m^{\l}) \in \Lnil$. 
Then $\Lnil$ is spanned by $U_{\e}^{+}v(m^{\l})$ 
as a $\bbC$-vector space. 
Therefore, by the similar way to the proof of $E_0$-case, 
we obtain $\widetilde{\rm{ev}}^{+}_{\bf{a}_{+}}(F_0)
=\widetilde{\rm{ev}}^{-}_{{\bf{a}}_{-}}(F_0)$ on $\Lnil$. 
\qed

\end{document}